\documentclass[11pt,a4paper]{preprint}
\usepackage[utf8]{inputenc}
\usepackage[english]{babel}
\usepackage{helvet, stmaryrd, ulem}
\usepackage{amsmath}
\usepackage{amsthm}
\usepackage{amsfonts}
\usepackage{amssymb}
\usepackage{mathtools}
\usepackage{mathrsfs}
\usepackage{dsfont}
\usepackage{tabularx}
\usepackage{xcolor}
\usepackage{bbm}
\usepackage{mhequ}
\usepackage{todonotes}
\usepackage{multirow}
\usepackage{caption}
\usepackage[paper=a4paper,left=25mm,right=25mm,top=25mm,bottom=35mm]{geometry}
\usepackage[noabbrev,capitalize]{cleveref}
\usepackage[nottoc]{tocbibind}
\usepackage{comment}

\newtheorem{theorem}{Theorem}[section]

\crefname{definition}{Definition}{Definitions}
\newtheorem{proposition}{Proposition}[section]
\crefname{proposition}{Proposition}{Propositions}
\newtheorem{lemma}{Lemma}[section]
\crefname{lemma}{Lemma}{Lemmata}

\crefname{corollary}{Corollary}{Corollaries}

\crefname{example}{Example}{Examples}
\newtheorem{assumption}{Assumption}[section]
\crefname{assumption}{Assumption}{Assumptions}

\theoremstyle{definition}
\newtheorem{remark}{Remark}[section]
\crefname{remark}{Remark}{Remarks}

\numberwithin{equation}{section}

\newcommand{\supp}{\operatorname{supp}}
\newcommand{\R}{\mathbb{R}}
\newcommand{\cR}{\mathcal{R}}

\newcommand{\T}{\mathbb{T}}
\newcommand{\N}{\mathbb{N}}
\newcommand{\E}{\mathds{E}}
\newcommand{\p}{\mathds{P}}

\newcommand{\F}{\mathcal{F}}
\newcommand{\cB}{\mathcal{B}}

\newcommand{\calC}{\mathscr{C}}

\renewcommand{\mathcal}[1]{\mathscr{#1}}
\renewcommand{\epsilon}{\varepsilon}
\newcommand{\eps}{\varepsilon}
\newcommand{\bF}{\mathbb{F}}
\newcommand{\Z}{\mathbb{Z}}
\newcommand{\C}{\mathbb{C}}

\newcommand{\squeeze}[2][0]{\mbox{$\medmuskip=#1mu\displaystyle#2$}}

\DeclarePairedDelimiter{\abs}{\lvert}{\rvert}
\DeclarePairedDelimiter{\norm}{\lVert}{\rVert}

\DeclarePairedDelimiter{\paren}{(}{)}

\begin{document}
\title{Higher order convergence rate in probability for approximations of the stochastic Burgers and Cahn-Hilliard equation}

	\date{\today}
\author{Ana Djurdjevac\footnote{University of Oxford, \email{ana.djurdjevac@maths.ox.ac.uk}}, Máté Gerencsér\footnote{TU Wien,
\email{mate.gerencser@asc.tuwien.ac.at}}, Helena Kremp\footnote{WIAS Berlin,
\email{kremp@wias-berlin.de}}}
\maketitle

	\begin{abstract}
We consider temporal approximations of stochastic PDEs driven by additive space-time white noise, with nonlinearities that may lose both integrability and regularity. Our key examples are the stochastic Burgers equation in dimension $1$ and the stochastic Cahn-Hilliard equation in dimensions $\leq 3$.
We show that the accelerated exponential Euler scheme achieves a rate of convergence in probability that is well beyond the state of the art: arbitrarily close to $3/4$ in the case of Burgers and $1-d/8$ in the case of Cahn-Hilliard.
		
		\bigskip
		
		 {{\sc Mathematics Subject Classification (2020):  60H15, 60H35 }
		
		}

	\end{abstract}

\section{Introduction}
 We consider numerical approximations of SPDEs posed on the torus $\mathbb{T}^d$ of the form
\begin{equ}\label{eq:SPDE}
    \partial_tu+(-\Delta)^{\sigma/2} u=\mathcal{D}_\alpha G(u)+\xi,
\end{equ}
where $\xi$ is space-time white noise on $\mathbb{T}^d$, $G:\R\to\R$ is a sufficiently regular function and $\mathcal{D}_\alpha$ is a differential operator of order $\alpha \in\N$, i.e. $\mathcal{D}_\alpha f=\sum_{|\ell|\leq\alpha}a_\ell\partial^\ell f$ with $a_{\ell}\in\R$ for multi-indices $\ell=(\ell_1,\dots,\ell_d)\in\N^d_0$ and $\partial^\ell=\partial^{\ell_1}_1\cdots\partial^{\ell_d}_d$. We assume that $\sigma>d$, so that the solution is a function in the space variable and therefore the nonlinearity $G(u)$ is well-defined.
Examples of interest are the stochastic Burgers equation, where
\begin{equ}\label{ex:Burgers}
    d=1,\quad\sigma=2,\quad\alpha=1,\quad G(u)=u^2,\quad \mathcal{D}_\alpha=\partial_x\tag{B}
\end{equ}
and the stochastic Cahn-Hilliard equation, where
\begin{equ}\label{ex:CH}
    d=1,2,3,\quad\sigma=4,\quad \alpha=2,\quad G(u)=u^3-u,\quad\mathcal{D}_\alpha=\Delta.\tag{CH}
\end{equ}

In the generality of \eqref{eq:SPDE}, obtaining strong convergence rates of full discretisations is rather challenging.
To our best knowledge, there is no proof for any scheme achieving positive rate in $L^p(\Omega)$ for \eqref{ex:Burgers} or for \eqref{ex:CH} in dimensions $d=2$ or $3$.
One difficulty is that the possibly superlinear growth of $G$ can cause certain temporal discretisation schemes to explode in $L^p(\Omega)$ (see \cite{beccari2019strong}).
However, since these blowups in $L^p(\Omega)$ are accumulated on events of very small probability, it is often possible to prove rates with another notion of strong convergence, namely convergence in probability.
For example, while there is no known temporal rate in $L^p(\Omega)$ for \eqref{ex:Burgers}, \cite{Printems} showed rate $1/4$ of convergence in probability of an exponential Euler scheme.

In the present paper we study the rate of convergence in probability of the accelerated exponential Euler scheme. This scheme was proposed in \cite{Jentzen-Kloden} in order to achieve higher temporal convergence rates. In the simplest case $\sigma=2$, $d=1$, $\alpha=0$, a proof for the numerically observed temporal rate \footnote{To avoid obscuring the overview of the literature with $\epsilon$-s, we do not distinguish between rate $\alpha$ and $\alpha-\epsilon$ for all $\epsilon>0$ in the mentioned results.}\label{fn} $1$ was recently established by the authors in \cite{DjGK} with respect to convergence in $L^{p}(\Omega)$. In the present paper we employ the strategy of \cite{DjGK} to equations of the more general form \eqref{eq:SPDE} in the context of convergence in probability as initiated in \cite{Printems}. The rates that we obtain substantially improve the best known temporal rates
(see \cref{table} below for a comparison).

\subsection{Literature overview}

We will first comment on possible strategies for improving the standard convergence rate in the strong sense, which, as mentioned in the introduction, cannot be extended in an obvious way to the general equation \eqref{eq:SPDE}, and in particular not to the cases \eqref{ex:Burgers}  and \eqref{ex:CH}.
The so-called accelerated exponential Euler scheme for SPDEs was introduced in \cite{Jentzen-Kloden},   to overcome rigorous lower bounds that apply to a large classes of strong approximations \cite{Davie-Gaines, BGJK}. In the context of \eqref{eq:SPDE} with $\sigma=2$, $d=1$, $\alpha=0$, in \cite{Jentzen-Kloden} the mapping $u\mapsto G\circ u$ was replaced by a mapping $u\mapsto F(u)$, $F:L^2(\T)\to L^2(\T)$, and the assumptions were formulated in terms of $F$. However, there exists no truely nonlinear function $G$ such that $F(u)=G\circ u$ satisfies the assumptions required in \cite{Jentzen-Kloden} (see e.g. \cite[Prop.~1.2,~Rem.~1.3]{DjGK}). It was later proved in \cite{Jentzen12} that the accelerated exponential Euler scheme does achieve a higher 
temporal convergence rate of $1/2$ in the case of globally Lipschitz nonlinearity $G$, and in \cite{Wang} the same temporal rate $1/2$ is proved for a tamed version of the scheme, in the case of the Allen-Cahn nonlinearity $G(u)=-u^3+u$. These rates improved the classical temporal rate of $1/4$ (see e.g. \cite{Gy, becker2023strong}), but did not match the empirically observed rate$^{\ref{fn}}$ $1$. The temporal rate $1$ was recently proved in \cite{DjGK}, both for the accelerated exponential Euler scheme (if $G$ has linear growth) and for a splitting version of the scheme (if $G$ has superlinear growth).
For other results on the accelerated exponential Euler scheme see  \cite{HutzenthalerJentzenSalimova2018, ChenCuiHongSheng2022,Brehier2024}.
For the stochastic Burgers equation \eqref{ex:Burgers}, existence and uniqueness of solutions goes back to \cite{DaPrato}.
Various spatial discretizations of the stochastic Burgers equation have been investigated, and a spatial convergence rate of order $1/2$ has been established. Initial results were obtained for finite difference schemes in \cite{Alabert2006}, with several later developments addressing both the choice of spatial norms and the discretization methods themselves, such as the spectral Galerkin approach in \cite{BJ13} and the Fourier spectral collocation method in \cite{MinooHosseini_zbMATH06107055}. It should be emphasized that the choice of finite difference scheme requires particular care: as demonstrated in \cite{HairerVoss_zbMATH06031653, HairerMaas12_zbMATH06067453} and further developed in later works, different discretizations may converge to different limiting equations.

Establishing convergence rates for temporal discretizations has proven to be considerably more challenging. Apart from the aforementioned work \cite{Printems} obtaining rate $1/4$ of convergence in probability, the works \cite{JSW19, MazzonettoSalimova20} proved convergence of fully discrete space–time schemes based on nonlinearity-truncated, accelerated exponential Euler methods, but without providing an explicit rate. Recently, in \cite{wang2024strong}, a fully explicit space–time scheme combining spectral Galerkin and tamed accelerated exponential Euler method was analyzed for the stochastic Burgers–Huxley equation driven by space–time white noise, yielding strong convergence rates. However, the additional cubic reaction term provides a dissipative effect that stabilizes the dynamics, which makes it easier to derive uniform moment bounds and to prove stability of explicit schemes.

Note furthermore that for equations driven by more regular noise the literature is much richer. In particular, for trace-class noise (see, e.g., \cite{HutzenthalerJentzenzbMATH07206753}) there are several results establishing both strong convergence rates \cite{hutzenthaler2019strong, hutzenthaler2022strong} and recent result on weak convergence rates \cite{brehier2024weak}, while  with fractional-type Gaussian noise with Hurst parameter  $H \in (1/2,1)$, strong convergence rates for exponential Euler and spectral Galerkin schemes were obtained in \cite{WangCao}.

For the stochastic Cahn-Hilliard equation \eqref{ex:CH}, existence and uniqueness of solutions goes back to \cite{DaPratoDebussche}. Higher regularity estimates were obtained in \cite{CW01}, that will be employed in our analysis and unlike in the Burger's case, do not immediately follow from semigroup estimates combined with the basic well-posedness result.

Several recent works have addressed the numerical approximation of the stochastic Cahn--Hilliard equation with space--time white noise.
In the one-dimensional setting, 
\cite{cai2023strong} show that a spectral Galerkin discretization in space coupled with a tamed exponential Euler scheme in time attains strong convergence rates of 
 order $3/8$ in time and $3/4$ in space.
For the same scheme, it was further shown in \cite{brehier2022weak} that the weak order of convergence is twice the corresponding strong order with respect to both spatial and temporal discretization parameters. 
A spectral Galerkin approximation in the one dimensional setting together with an accelerated implicit Euler scheme in time was studied in \cite{cui2020}, 
yielding a strong convergence rate of order $3/4$ both in space and time. A semi-discrete method based on an implicit Euler time discretization combined with a consistent noise approximation, 
and its convergence to the sharp-interface limit, was investigated in \cite{BanasMukam}. 
A fully discrete scheme for fourth-order SPDEs (which corresponds to $\sigma=4$ and $\alpha=0$ in our setting) with a bounded and regular nonlinearity on the $d$-dimensional torus for $d=1,2,3$, combining an accelerated Euler scheme in time with a spectral Galerkin discretization, was analyzed in \cite{blomker2025numerical}, 
where spatial and temporal convergence rates of order $2-\frac{d}{2}$ and $\frac{3}{2}-\frac{d}{4}$, respectively, were obtained using the sewing techniques\footnote{We remark however, that at least in dimension $d=1$ the claimed rate is not proved, as the regularity bound \cite[Lemma~3.4]{blomker2025numerical} is false for $\lambda>2$.}.
For the one-dimensional stochastic Cahn--Hilliard equation driven by
additive space--time white noise, \cite{cui2020} obtain an almost
$3/4$-order strong temporal convergence rate for an accelerated
implicit Euler scheme, whereas the backward Euler scheme of \cite{cai2023strong} and the explicit tamed exponential Euler scheme of \cite{QiCaiWang2022} converge with almost order $3/8$.
For the space--time white noise in dimensions $d=2,3$, these
results are not applicable, and no directly comparable strong temporal
rate for the classical nonlinear equation appears to be available.

\begin{center}
    
\begin{tabular}{|c|c|c|c|}
 \hline
 & $d$ & State of the art & \cref{thm:main} \\
\hline
Burgers    & 1 & 1/4   \cite{Printems}           & 3/4          \\
\hline
\multirow{3}{*}{Cahn-Hilliard} & 1 & 3/4 \cite{cui2020}             & 7/8          \\
                               & 2 & N/A             & 3/4          \\
                               & 3 & N/A              & 5/8  \\\hline       
\end{tabular}
\captionof{table}{Comparison of previous best and the new temporal rates}\label{table}
\end{center}

In this paper, we focus solely on the temporal discretization and convergence rates.
Concerning the spatial rate, for the stochastic Burgers equation, \cite{BJ13, Alabert2006} prove an almost $1/2$-order spatial convergence rate, uniformly in time. \cite{Alabert2006} obtain this rate for a finite-difference approximation in $L^2$, while \cite{BJ13} establish it for a spectral Galerkin approximation in the stronger $L^\infty$-norm. \cite{HairerMatetski_zbMATH06600872} further show that, for a broad class of spatial approximation schemes, the optimal rate is almost $1/2$ in the uniform norm, although the limiting equation may contain an additional scheme-dependent correction term. For the stochastic Cahn--Hilliard equation with additive space--time white noise in one dimension, \cite{cui2020,cai2023strong} prove strong spatial convergence of order $3/2$.

\subsection{Proof strategy and power counting for guessing the temporal rate}
Let $(P_{t}^{\sigma})_{t\geq 0}$ denote the semigroup of $-(-\Delta)^{\sigma/2}$ with $P^{\sigma}_t f= p_{t}^{\sigma}\ast f$ (see below for the precise definition of the semigroup) and let
\begin{equ}\label{eq:OU}
    O_t^{\sigma}(x) := \int_{0}^{t} \int_{\mathbb{T}^d} p_{t-s}^{\sigma}(x-y)\xi(ds,dy).
\end{equ}
By classical regularity bounds for the stochastic convolution, we have that for any $\theta\in [0,\frac{\sigma}{2})$, almost surely $O^{\sigma}\in C_{T}^{\frac{\theta}{\sigma}}\calC^{\frac{\sigma-d}{2}-\theta}$ (cf. \cref{prop:ou} below).
The analysis of the temporal rate reduces to the estimation of the following quadrature term:
\begin{align}\label{eq:quadrature}
    \int_{0}^{t} P_{t-s}^{\sigma}\paren[\Big]{f_{R}(O_{k_{M}(s)}^\sigma)-f_R(O_{s}^\sigma)}\,ds
\end{align}
in $L^{p}(\Omega,L^{\infty}(\T^d))$, where $f_R$ is an appropriate localization of the nonlinearity $f$ and $k_{M}(s)$ denotes the largest time grid point less than or equal to $s$.
A first approach to estimate $\eqref{eq:quadrature}$ consists of writing
\begin{align*}
    f_R(O_{k_{M}(s)}^\sigma)-f_R(O_{s}^\sigma) = \mathcal{D}_{\alpha}\paren[\Bigg]{\paren[\Big]{\int_0^1 \partial G_R ( \lambda O_{k_{M}(s)}^\sigma + (1-\lambda) O_{s}^\sigma ) d\lambda  } (O_{k_{M}(s)}^\sigma-O_{s}^\sigma)}.
\end{align*}
By product estimates in Besov spaces, the fact that $\mathcal{D}_{\alpha}$ decreases regularity by $\alpha$, as well as the regularity estimates for $O^{\sigma}$ we deduce that for small $\epsilon>0$, under the assumption that $d< 2\alpha$, $\sigma\geq 2\alpha$ and $\sigma < 2\alpha+d$,
\begin{align*}
    \norm{f_R(O_{k_{M}(s)}^\sigma)-f_R(O_{s}^\sigma)}_{\calC^{-\frac{\sigma}{2}+\epsilon}}\lesssim \norm{O^{\sigma}}_{C_T\calC^{\frac{\sigma}{2}-\alpha}}\norm{O^{\sigma}}_{C_T^{1-\frac{d+2\alpha}{2\sigma}-\frac{\epsilon}{\sigma}}\calC^{-\frac{\sigma}{2}+\alpha+\epsilon}} (s-k_{M}(s))^{1-\frac{d+2\alpha}{2\sigma}-\frac{\epsilon}{\sigma}}.
\end{align*}
Next, we employ the regularisation of the semigroup $P^{\sigma}$ in order to improve the latter regularity $-\frac{\sigma}{2}+\epsilon$ to $\epsilon$ and end up estimating $\eqref{eq:quadrature}$ in $L^{\infty}(\T^d)$ by $t^{\frac{1}{2}+\frac{\epsilon}{\sigma}} M^{-1+\frac{d+2\alpha}{2\sigma}+\frac{\epsilon}{d}}$.
We propose an alternative estimate for \eqref{eq:quadrature} using a stochastic sewing strategy, which allows to trade time regularity for an increased temporal rate, such that we obtain a bound by $t^{\frac{\epsilon}{\sigma}} M^{-\frac{3}{2}+\frac{d+2\alpha}{2\sigma}+\frac{\epsilon}{\sigma}}$.

In the below we detail our heuristics for guessing that particular temporal rate following the methods in \cite{DjGK}. 
The following heuristics require familiarity with stochastic sewing with blow-up/weights, see Lemma \ref{lem:ss} below. Readers who are not familiar with these techniques may safely omit the below paragraph.
\begin{itemize}
    \item Recall that $O^{\sigma}\in C_{T}^{\frac{\theta}{\sigma}}\calC^{\frac{\sigma-d}{2}-\theta}$ for any $\theta\in [0,\frac{\sigma}{2})$. The solution $u$ inherits the same regularity and we assume that $\sigma>d$, so that the solution is a function. In particular, we have that $u\in C_{T}^{\frac{1}{2}-\frac{d}{2\sigma}}\calC^{0}$.
    \item The nonlinearity thus satisfies $f_R(u)=\mathcal{D}_{\alpha} G_R(u)\in C_{T}^{\frac{1}{2}-\frac{d}{2\sigma}+\frac{\lambda}{\sigma}}\calC^{-\alpha-\lambda}$ for $\lambda\in [0,\frac{d}{2})$.
    \item We measure the error between $u_t$ and its numerical scheme in $L^{\infty}(\T^d)$. 
    Next, we can use the regularisation coming from the semigroup to upgrade the spatial regularity of the nonlinearity, that is $\calC^{-\alpha-\lambda}$, to $L^{\infty}$, at the expense of a temporal blow-up. It turns out, that this blow-up is not allowed to exceed $1/2$ due to restrictions in the sewing lemma. This gives the relation
    \begin{equ}\label{eq:lambda}
        \frac{\alpha+\lambda}{\sigma}=\frac{1}{2}\qquad\Leftrightarrow\qquad \lambda=-\alpha+\frac{\sigma}{2}.
    \end{equ}
   Recalling the restriction 
   $\lambda\in [0,d/2)$, we obtain a relation of $\alpha, \sigma$ and $d$ of $\sigma<d+2\alpha$ and $\sigma\geq 2\alpha$.
  
    \item The temporal rate is given by the time regularity of the nonlinearity for the choice of $\lambda$ from \eqref{eq:lambda} plus $\frac{1}{2}$ that we gain due to the application of the stochastic sewing lemma:
    \begin{equ}
        \frac{1}{2}+\frac{1}{2}-\frac{d}{2\sigma}+\frac{\lambda}{\sigma}=1-\frac{d}{2\sigma}-\frac{\alpha}{\sigma}+\frac{1}{2}=\frac{3}{2}-\frac{d+2\alpha}{2\sigma}\leq 1,
    \end{equ}
    where the rate is always bounded by $1$ due to the upper bound on $\lambda$.
    Thus for \eqref{ex:Burgers} we obtain a temporal rate $\frac{3}{4}$ and for \eqref{ex:CH} we obtain a temporal rate $1-\frac{d}{8}$ for $d=1,2,3$.
\end{itemize}

\begin{remark}
The stochastic Cahn-Hilliard equation, that Cook in \cite{Cook1970} originally provided the physical derivation for, is in fact driven by the spatial derivative\footnote{To be precise, in higher dimension one should take $(-\Delta)^{1/2}\xi$.} of a space-time white noise. Such an equation is only solvable classically for $d=1$, for the rigorous mathematical analysis see \cite{ElezovicMikelic1991, DebusscheZambotti2007}.
For $d\geq2$ one expects that  renormalisation is necessary, see e.g. \cite{BanasYangZhu2023}.

This example could motivate a further generalisation, allowing a noise of the form $(-\Delta)^\gamma\xi$ in place of $\xi$ in \eqref{eq:SPDE}.
The $1$-dimensional physical Cahn-Hilliard-Cook  equation \cite{Cook1970} would then correspond to
\begin{equation}\label{ex:CHC}
    d=1,\quad\sigma=4,\quad \alpha=2,\quad\gamma=\frac{1}{2},\quad G(u)=u^3-u,\quad\mathcal{D}_\alpha=\Delta.\tag{CHC}
\end{equation}
Although it is not expected that this changes anything essential in the arguments of the present paper, we refrain from this generalisation for two reasons. First, carrying another parameter $\gamma$ would make the already at times lengthy expressions in the exponents even longer.
Second, the a priori regularity that we assume for the equations in \cref{ass:1} does not easily follow from the existing literature for \eqref{ex:CHC},
although remarks in \cite{DebusscheZambotti2007} indicate that the methods therein should be adaptable.
One can nevertheless perform the power counting heuristics above and make the educated guess that the methods of the present paper should give temporal rate $5/8$ for \eqref{ex:CHC}.
\end{remark}

The paper is organised as follows. In \cref{sec:setting} we state our main result, \cref{thm:main}. In \cref{sec:prelim} we summarize all needed heat-kernel and functional estimates, state a priori estimates for the stochastic convolution, as well as prove auxiliary lemmata. Furthermore \cref{sec:prelim} features a priori estimates for our temporal scheme and the solution in the case that the nonlinearity is localized, which are the only a priori estimates that are needed for our methods. In \cref{sec:proof} we prove our main \cref{thm:main}, while in \cref{sec:ex} we show that our assumptions are satisfied in the model cases of Burgers and Cahn-Hilliard equation.

\bigskip

\section*{Acknowledgments}
The authors thank Dirk Blömker for discussions related to the problem. 
\begin{itemize}
    \item ADj: ADj gratefully acknowledges  funding by the Deutsche Forschungsgemeinschaft (DFG, German Research Foundation) CRC/TRR 388 "Rough Analysis, Stochastic Dynamics and Related Fields“ – Project ID 516748464 and by Daimler and Benz Foundation as part of the scholarship program for junior professors and postdoctoral researchers. 
    \item MG: Funded by the European Union (ERC, SPDE, 101117125). Views and opinions expressed are however those of the author(s) only and do not necessarily reflect those of the European Union or the European Research Council Executive Agency. Neither the European Union nor the granting authority can be held responsible for them.
    \item HK: This research was funded in whole or in part by the Austrian Science Fund (FWF) [10.55776/P34992]. HK acknowledges support from the DFG CRC/TRR 388 ``Rough Analysis, Stochastic Dynamics and Related Fields" Project ID 516748464 within subproject A02.
\end{itemize}
\section*{AI Declaration}
GPT-6 was solely used for a final proofreading of the manuscript in order to rule out typos or minor inconsistencies.

\section{Setting and statement}\label{sec:setting}
Let $(\Omega,\F,\p)$ be a probability space. We fix a time horizon $T>0$ and a dimension $d\in\N$. The $d$-dimensional torus $\T^d$ is defined as $\T^d=\R^d/\Z^d$.
The space-time white noise $\xi$ is defined as a mapping from the Borel sets $\cB([0,T]\times\mathbb{T}^d)$ into $L^{2}(\Omega)$, such that for any collection $A_{1},\dots,A_{k}\in \cB([0,T]\times\mathbb{T}^d)$, the vector $(\xi(A_{1}),\dots,\xi(A_{k}))$ is Gaussian with zero mean and covariance $\E[\xi(A_{i})\xi(A_{j})]=\lambda(A_{i}\cap A_{j})$, where $\lambda$ denotes the Lebesgue measure.
We also fix a filtration $\bF=(\F_{t})_{t\in[0,T]}$, such that $(\Omega,\F,\bF,\p)$ is complete and such that for any $t\in[0,T]$, $A\in \cB([0,t]\times\mathbb{T}^d)$, $B\in \cB([t,T]\times\mathbb{T}^d)$, $\xi(A)$ is $\F_{t}$-measurable and $\xi(B)$ is independent of $\F_{t}$. An example for $\bF$ would be the completed filtration generated by $\xi$. 
By stochastic integral of (deterministic) square integrable functions against $\xi$ we mean the unique isometric and linear extension of the map $\mathbf{1}_A\mapsto\xi(A)$ to $L^2([0,T]\times\T^d)$.

For $k\in\Z^d$, denote the Fourier modes on $\T^d$ by $e_{k}(x)=e^{2\pi ik\cdot x}$. The set $(e_k)_{k\in\Z^d}$ is an orthonormal basis of $L^2(\T^d,\C)$.
For $f\in L^1(\T^d,\C)$, its Fourier transform is denoted by $\F f (k)=\hat{f}(k)=\int_{\mathbb{T}^d}e^{-2\pi i k\cdot x}f(x)dx$, $k\in\Z^d$. 
Denote by $(P_{t}^{\sigma})_{t\geq 0}$ the semigroup generated by $-(-\Delta)^{\sigma/2}$ on the torus, that is
$$P_{t}^{\sigma}f=\F^{-1}(e^{-\abs{2\pi k}^{\sigma}t}\F f (k)).$$
One can equivalently write $P_t^\sigma f=p_t^\sigma\ast f$, where
$$
p_{t}^\sigma(x)=\sum_{k\in\mathbb{Z}^d}e^{-\abs{2\pi k}^{\sigma}t}e^{2\pi i k \cdot x}.
$$ 
\noindent We consider the mild formulation of \eqref{eq:SPDE}: 
\begin{align}\label{eq:mild}
    u_{t}=P_{t}^{\sigma}u_{0}+\int_{0}^{t}P_{t-s}^{\sigma}\big(\mathcal{D}_{\alpha}G(u_{s})\big)ds+O^\sigma_t,
\end{align}
where the stochastic convolution $O^\sigma$, also referred to as the Ornstein-Uhlenbeck (OU) process, is defined in \eqref{eq:OU}.
Let here and below $\lceil r\rceil$ denote the smallest natural number $n\geq r$.

Throughout the whole paper we assume the following two assumptions.
\begin{assumption}\label{ass:G}
    We assume that $G\in C^{\lceil \frac{\sigma-d}{2}\rceil + 2}(\R)$, i.e. $G:\R\to \R$ is $m$-times continuously differentiable for $m=\lceil \frac{\sigma-d}{2}\rceil + 2$.
\end{assumption}

\begin{assumption} \label{ass:1}
We assume that $\sigma,\alpha, d$ are such that $\sigma\geq 2\alpha$ and $d<\sigma <d+2\alpha$ and $2\alpha>d$. Moreover we assume that for any $\epsilon\in (0,1)$, $\norm{u_0}_{\calC^{\frac{\sigma-d}{2}-\epsilon}(\T^d)}\lesssim \mathcal{M}(\epsilon)<\infty$ and that there exists a unique mild solution $u$ to \eqref{eq:SPDE} which satisfies, that for any $\lambda\in [0,\frac{\sigma}{2})$ and any small enough $\epsilon>0$,  $u\in C_{T}^{\lambda/\sigma}\calC^{\frac{\sigma-d}{2}-\lambda-\epsilon}(\mathbb{T}^d)$ almost surely.
In particular, $u$ is a function in time-space $(t,x)$. 
\end{assumption}
\begin{remark}
    We note that both assumptions on the exponents are satisfied if $\sigma=2\alpha >d$ and in particular for \eqref{ex:Burgers} and \eqref{ex:CH}. Moreover the existence and uniqueness of a local solution with the prescribed regularity is automatic (cf. \cite[Theorem 7.8]{hairer2023}). Since $G$ is allowed to grow superlinearly, the fact that the solution is global is not automatic and
    will be verified case by case, see \cref{prop:wp} below.
\end{remark}
\begin{remark}
   The case of $\alpha=0$ is excluded in \cref{ass:1}. We refer to \cite{DjGK} for the case $\alpha=0$, $\sigma=2$, $d=1$, in which a temporal rate of (almost) $1$ was proven (with convergence in $L^p(\Omega)$ for bounded nonlinearities $f$). For $\sigma>2$ the stochastic convolution is even more regular and the arguments of \cite{DjGK} for $d=1$ apply to show a temporal rate of (almost) $1$. For $\alpha=0$, $\sigma\in (1,2)$, $d=1$, we expect a temporal rate of (almost) $\frac{3}{2}-\frac{1}{\sigma}$ (without proof, but one may adapt the arguments of \cite{DjGK} to show the latter). Indeed, in the power counting heuristics above, whenever $\alpha=0$ and $2d\geq \sigma$, one chooses the blow up $\lambda=\frac{\sigma-d}{2}$, which results in a temporal rate $\frac{3}{2}-\frac{d}{\sigma}$.
\end{remark}

\bigskip

As a temporal discretization, we consider the accelerated exponential Euler scheme as follows. We define the scheme $V^M$ iteratively for the time grid points $t_k=\frac{Tk}{M}=k h$ for $k=0,\dots M$ by letting $V^{M}_{0}=u_{0}$ and
\begin{align*}
    V^{M}_{t_{k+1}}=P_{h}^{\sigma}V^{M}_{t_{k}}  -(-\Delta)^{-\sigma/2}(P_h^{\sigma}-\operatorname{Id})\mathcal{D}_{\alpha}G(V^{M}_{t_k}) + O_{t_{k+1}}^{\sigma}-P_{h}^{\sigma}O_{t_k}^{\sigma}, \quad k=0,\dots, M-1.
\end{align*}
Note that $-(-\Delta)^{\sigma/2}$ is invertible only on $\mathrm{span}\{e_j:\,j\in \mathbb{Z}\setminus\{0\}\}$, therefore we set by convention  $-( -\Delta)^{-\sigma/2}(P^\sigma_h-\mathrm{Id})e_0\coloneqq h e_0$.
On grid points the inductive form agrees with the following mild form, which we further take as a definition for any continuous time point $t\in [0,T]$, with the notation $k_{M}(t):=\lfloor h^{-1}t\rfloor h$,
\begin{align}\label{eq:scheme1}
    V^{M}_{t}=P_{t}^{\sigma}u_{0}+\int_{0}^{t}P_{t-s}^{\sigma}\mathcal{D}_{\alpha}G(V^{M}_{k_{M}(s)})ds+O_t^{\sigma},\quad t\in [0,T].
\end{align}
We can finally state our main theorem.
\begin{theorem}\label{thm:main}
Let $T>0$ and let \cref{ass:G} and \cref{ass:1} hold. Let $u$ be the mild solution of \eqref{eq:SPDE}.
Let $V^{M}$ be the approximation obtained by the accelerated exponential Euler scheme defined in \eqref{eq:scheme1}.Then for any $\epsilon>0$, we have convergence of the scheme in probability with rate $M^{-\frac{3}{2}+\frac{d+2\alpha}{2\sigma}+\epsilon}$, i.e. for any $\delta>0$,
\begin{align}\label{eq:Linfty}
    \limsup_{M\to\infty}\p\Big(\sup_{t\in[0,T]}\norm{u_t-V^{M}_t}_{L^{\infty}}\geq \delta M^{-\frac{3}{2}+\frac{d+2\alpha}{2\sigma}+\epsilon}\Big)=0.
\end{align}
\end{theorem}
\begin{remark}
A similar statement holds for the full discretisation. Indeed, define the time-space discretised scheme by 
$V^{M,N}$ with $$V^{M,N}_t=P^{\sigma,N}_t u_0 +\int_{0}^{t}P^{\sigma,N}_{t-s} \mathcal{D}_{\alpha}G(V^{M,N}_{k_{M}(s)})ds + O^{\sigma,N}_t, \quad t\in[0,T],$$ where $P^{\sigma,N}:=\Pi_N P^{\sigma}=P^{\sigma}\Pi_N$ and $O^{\sigma,N}=\Pi_N O^{\sigma}$ for the projection $\Pi_N$ onto the first $N$ Fourier modes. Let furthermore $u^N$ be the solution to the Galerkin projected equation, that is for $P^{\sigma}$ replaced by $P^{\sigma,N}$ and $O^{\sigma}$ replaced by $O^{\sigma,N}$. Then we have that for any $\delta>0$,
\begin{align}\label{eq:Linfty-G}
 \limsup_{M\to\infty}\sup_{N\in\N}\p\Big(\sup_{t\in[0,T]}\norm{u_t^N-V^{M,N}_t}_{L^{\infty}}\geq \delta M^{-\frac{3}{2}+\frac{d+2\alpha}{2\sigma}+\epsilon}\Big)=0.
\end{align}
The proof is completely analogue to the proof of \cref{thm:main}, using the semigroup estimates for $(P_t^{\sigma,N})$ instead of $(P_t^{\sigma})$ and regularity bounds for $O^{\sigma,N}$ instead of $O^{\sigma}$, uniformly in $N$.
\end{remark}

\section{Preliminaries}\label{sec:prelim}
We introduce the Besov spaces as follows. Let $(\rho_{j})_{j\geq -1}$ be a smooth dyadic partition of unity, i.e. a family of functions $\rho_{j}\in C^{\infty}_{c}(\R^d)$ for $j\geq -1$, such that 
\begin{itemize}
\item $\rho_{-1}$ and $\rho_{0}$ are non-negative even functions such that the support of $\rho_{-1}$ is contained in $B_{1/2}$, the ball of radius $1/2$ around $0$, and the support of $\rho_{0}$ is contained in $B_1\setminus B_{1/4}$;
\item $\rho_{j}(x)=\rho_{0}(2^{-j}x)$, $x\in\R^d$, $j\geq 0$;
\item $\sum_{j=-1}^{\infty}\rho_{j}(x)=1$ for every $x\in\R^d$;
\item $\operatorname{supp}(\rho_{i})\cap \operatorname{supp}(\rho_{j})=\emptyset$ for all $\abs{i-j}>1$.
\end{itemize}
The existence of such a partition of unity is classical (see e.g. \cite{Bahouri2011}).
We denote by $\mathcal{S}'(\mathbb{T}^d)$ the space of Schwartz distributions on the torus (i.e. the dual of  $\mathcal{S}(\mathbb{T}^d):=C^\infty(\T^d)$).
Note that for any $u\in \mathcal{S}'(\T^d)$ its Fourier transform is meaningful, and therefore one can define the operators (also known as Littlewood-Paley blocks), $j\geq -1$, 
\begin{align}
    \Delta_j:\mathcal{S}'(\T^d)\to C^\infty(\T^d)\qquad \Delta_j u=\mathcal{F}^{-1}\big(k\mapsto \rho_{j}(k)\mathcal{F}(u)(k)\big).
\end{align}
We then define the Besov spaces on the torus $\T^d$ for $p,q\in [1,\infty]$, $\theta\in\R$
\begin{align}\label{def:bs}
B^{\theta}_{p,q}:=\big\{ u\in\mathcal{S}'(\mathbb{T}^d):\norm{u}_{B^{\theta}_{p,q}}=\norm{(2^{j\theta}\norm{\Delta_{j}u}_{L^{p}(\T^d)})_{j\geq -1}}_{\ell^{q}}<\infty\big\}.
\end{align}
We introduce the shorthand $\calC^{\theta}:=B^\theta_{\infty,\infty}$ for $\theta\in\R$.
We collect properties of Besov spaces below, which are relevant to us. 
\begin{itemize}
    \item (H\"older spaces) If $\theta\in\R^{+}\setminus\N$, then $\calC^{\theta}$ coincides with the space of functions whose derivatives up to order $\abs{\beta}\leq\lfloor \theta\rfloor$ are bounded and whose derivatives of order $\lfloor \theta\rfloor$ are $(\theta-\lfloor \theta\rfloor)$-H\"older continuous and one has the norm equivalence
    \begin{align*}
       \norm{u}_{\calC^{\theta}}\simeq \sum_{\abs{\beta}\leq\lfloor \theta\rfloor}\norm{\partial^{\beta} u}_{\infty}+ \sum_{\abs{\beta}=\lfloor \theta\rfloor} \sup_{x\neq y}\frac{\abs{\partial^{\beta}u(x)-\partial^{\beta}u(y)}}{\abs{x-y}^{\theta-\lfloor \theta\rfloor}}
    \end{align*}
    (cf.~\cite[Sec. 2.7, Examples]{Bahouri2011}), where $\beta=(\beta_1,\dots,\beta_d)\in\N^d$ with $\abs{\beta}=\beta_1+\dots+\beta_d$ and $\partial^\beta=\partial^{\beta_1}_1\cdots\partial^{\beta_d}_d$.
    \item (Embedding) If $1\leq p_{1}\leq p_{2}\leq \infty$, $1\leq q_{1}\leq q_{2}\leq\infty$ and $\theta\in\R$, then $B^{\theta}_{p_{1},q_{1}}$ is continuously embedded in $B_{p_{2},q_{2}}^{\theta-d(1/p_{1}-1/p_{2})}$ and we write $B^{\theta}_{p_{1},q_{1}}\hookrightarrow B_{p_{2},q_{2}}^{\theta-d(1/p_{1}-1/p_{2})}$ (cf. \cite[Proposition 2.71]{Bahouri2011}).
    \item (Derivatives) If $u\in B^{\theta}_{p,q}$ and $\beta\in\N^d_0$ is a multi-index, there exists a constant $C=C(\theta,\beta)$ such that
    \begin{align}\label{eq:derivative}
    \norm{\partial^{\beta}u}_{B^{\theta-\abs{\beta}}_{p,q}}\leq C\norm{u}_{B^{\theta}_{p,q}},
    \end{align} which follows from Bernstein inequality \cite[Lemma 2.1]{Bahouri2011}.
    \item (Products) If $\theta,\vartheta$ are such that  $\theta + \vartheta > 0$ and $p_1,p_2,q, p\in[1,\infty]$ with $\frac{1}{p}=\frac{1}{p_1}+\frac{1}{p_2}$, then for any two distributions $u\in B^{\theta}_{p_1,q}$ and $v\in B^{\vartheta}_{p_2,q}$ their product $uv$ is well-defined and
there exists a constant  $C=C(\theta,\vartheta,p_1,p_2,q)$, such that the bound
\begin{align}\label{eq:product}
    \Vert uv \Vert_{B^{\min(\theta,\vartheta)}_{p,q}}\leq C \Vert u\Vert_{B^{\theta}_{p_1,q}}\Vert v\Vert_{B^{\vartheta}_{p_2,q}}
\end{align} 
holds    (cf.~\cite[Section 2]{Bahouri2011}).
\item (Fractional heat kernel bounds) If $\theta\geq 0$, $\tilde\theta\in (0,\sigma]$, $\sigma>0$, then there exist constants $C_1(\theta)$ and $C_2(\tilde \theta)$ such that for any $\beta\in \R$, $p\in [1,\infty]$, $u\in B^\beta_{p,\infty}$,  and $t\in(0,1]$, the bounds
\begin{align}\label{eq:HK-1}
     \norm{P_{t}^{\sigma}u}_{B^{\beta+\theta}_{p,\infty}}&\leq C_1 t^{-\theta/\sigma}\norm{u}_{B^{\beta}_{p,\infty}},\quad \norm{(\operatorname{Id}-P_{t}^{\sigma})u}_{B^{\beta-\tilde\theta}_{p,\infty}}\leq C_2 t^{\tilde \theta/\sigma}\norm{u}_{B^{\beta}_{p,\infty}}
\end{align}
hold (cf. \cite[Lemma A.5, A.7]{GIP} whose proof adapts to $p\in[1,\infty]$ similarly as for $p=\infty$ and any $\sigma>0$ is covered, since $\F (e^{-\abs{z}^{\sigma}})\in L^1$, as well as \cite[Lemma 2.1]{Bahouri2011} together with $(\operatorname{Id}-P_{t}^{\sigma}) u=\int_{0}^{t} (-\Delta)^{\sigma/2} P_{s}^{\sigma} u ds $ from which the bound for $\operatorname{Id}-P_{t}^{\sigma}$ follows). 
\end{itemize}

For function spaces on $\R$, we only need the simple notion of $C^k_b$ for $k=0,1,\ldots$, denoting the space of bounded measurable functions whose distributional derivative up to order $k$ are essentially bounded, equipped with the canonical norm (note in particular that elements of $C^0_b$ are not assumed to be continuous).
Moreover denote by $(P_{t}^{\mathbb{R}})_{t\geq 0}$ the heat semigroup on $\mathbb{R}$, that is 
\begin{align}
    P_{t}^{\mathbb{R}}f=p_{t}^{\R}\ast f,\qquad\qquad p_{t}^{\R}(x)=\frac{1}{\sqrt{2\pi t}}e^{-x^2/2t}.
\end{align}
The following estimate is straight-forward: for any $0\leq s\leq t\leq 1$ one has
\begin{align}\label{eq:HK}
    \norm{(P_{t}^\R-P_s^\R)u}_{C^0_b}\leq |t-s|\norm{u}_{C^2_b}.
\end{align}
For a Banach space $X$, $C_{T}X$ denotes the space of continuous functions in time with values in $X$ equipped with the supremum norm. For $\gamma\in (0,1]$, define $$C_{T}^{\gamma}X : =\big\{u\in C_TX:\, \norm{u}_{C_{T}^{\gamma}X}:=\sup_{t\in[0,T]}\norm{u_{t}}_{X}+\sup_{0\leq s<t\leq T}\frac{\norm{u_{t}-u_{s}}_{X}}{(t-s)^{\gamma}}<\infty\big\}.$$
Below we use the notation $a\lesssim b$ if there exists a constant $C>0$, such that $a\leq Cb$.
The dependence of the constant $C$ will be clear from the context.
If we want to stress dependence of the constant $C$ on a parameter $\tau$, we write $a\lesssim_{\tau}b$.
Next, we collect some auxiliary results used in the sequel.

\begin{proposition}\label{prop:ou}
   Let $(O_{t}^{\sigma})_{t\geq 0}$ be as in \eqref{eq:OU}. 
     Let $p\in[1,\infty)$, $\lambda\in [0,\frac{\sigma}{2})$ and $\epsilon\in (0,\frac{1}{2})$. Then there exists a constant $C=C(p,T,\lambda,\epsilon)$, such that for all $0\leq s\leq t\leq T$ one has
    \begin{align}\label{eq:ou-moment-bounds}
        \E\norm{O_{t}^{\sigma}-O_{s}^{\sigma}}_{\calC^{\frac{\sigma-d}{2}-\lambda-\epsilon}}^p\leq C (t-s)^{\frac{p\lambda}{\sigma}}.
    \end{align}
    Furthermore, one has
    \begin{align}\label{eq:ou-reg}
        \E\norm{O^{\sigma}}_{C_{T}^{\frac{\lambda}{\sigma}}\calC^{\frac{\sigma-d}{2}-\lambda-\epsilon}}^{p}\leq C.
    \end{align}
Moreover if $\lambda\in [0,\frac{\sigma}{2})$, $\epsilon\in (0,\frac{1}{2})$, $\theta\in[0,\sigma-\lambda]$, then there exists a constant $C=C(p,T,\theta,\epsilon,\lambda)$ such that  for all $s,t\in[0,T]$ one has
    \begin{align}\label{eq:s-t-bound}
        \paren{\E\sup_{r\in[0,T]}\norm{(P_{t+s}^{\sigma}-P_{t}^{\sigma})O_{r}^{\sigma}}_{\calC^{\frac{\sigma-d}{2}-\lambda-\epsilon}}^{p}}^{1/p}\leq C t^{-\frac{\theta}{\sigma}}s^{\frac{\lambda+\theta}{\sigma}}.
    \end{align}
\end{proposition}
\begin{proof}
    The bounds \eqref{eq:ou-moment-bounds} and \eqref{eq:ou-reg} are classical. The argument is carried out in \cite[Proposition 3.1]{DjGK} for $\sigma=2$ and can easily be adapted to the case of general $\sigma> 0$. \eqref{eq:s-t-bound} follows from \eqref{eq:ou-reg} and the semigroup estimates with $P_{t+s}^{\sigma}-P_{t}^{\sigma}=P_{t}^{\sigma}(P_{s}^{\sigma}-\operatorname{Id})$.
\end{proof}

The next lemma collects useful properties on the nonlinearity $f(u)=\mathcal{D}_\alpha G(u)$. 
In order to obtain global Lipschitz estimates, we truncate $f$. To that aim, we define, for a compactly supported $\varphi\in C^{\infty}_{c}(\R)$ with $\supp(\varphi)\subset (-2,2)$ and $\varphi(x)=1$ for $\abs{x}\leq 1$ and a parameter $R>0$ with $\varphi_{R}(x)=\varphi(x/R)$, $x\in\R$, the truncated nonlinearities
\begin{align}\label{eq:trunc-f}
G_{R}(x):=G(\varphi_{R}(x)x),\quad x\in\R,\quad f_{R}(u):=\mathcal{D}_{\alpha} G_{R} (u), \quad u\in L^{\infty}.
\end{align}
\begin{lemma}\label{lem:fR} 
Let for $R>0$, $f_R$ be as in \eqref{eq:trunc-f}. 
We have that $f_R : L^{\infty}\to\calC^{-\alpha}$ is globally Lipschitz continuous, that is, there exists a constant $C_{R}> 0$, such that for all $u,v\in L^{\infty}$,
    \begin{align}\label{eq:fR-Lip}
        \norm{f_{R}(u)-f_{R}(v)}_{\calC^{-\alpha}}\leq C_{R}\norm{u-v}_{L^{\infty}}
    \end{align}
and bounded in the sense that for all $u\in L^{\infty}$,
\begin{align}\label{eq:fR-bdd}
        \norm{f_R(u)}_{\calC^{-\alpha}}\leq C_R.
\end{align}
Let $s\in (0, \lceil\frac{\sigma-d}{2}\rceil)\setminus\N$. 
Then $f_R:\calC^{s}\to\calC^{s-\alpha}$ is locally Lipschitz continuous, that is, for all $u,v\in \calC^{s}$,
\begin{align}\label{eq:fR-lLip}
        \norm{f_{R}(u)-f_{R}(v)}_{\calC^{s-\alpha}}\leq C_{R}(1+\norm{u}_{\calC^{s}}^{\lceil s\rceil }+\norm{v}_{\calC^{s}}^{\lceil s\rceil })\norm{u-v}_{\calC^{s}}
    \end{align}
    and  admits a polynomial growth, that is, for all $u\in \calC^{s}$,
\begin{align}\label{eq:fR-growth}
        \norm{f_R(u)}_{\calC^{s-\alpha}}\leq C_R(1+\norm{u}_{\calC^{s}}^{\lceil s\rceil })
\end{align}
for a (possibly different) constant $C_R>0$.
Furthermore, for any $\epsilon\in (0,2s)$ and for all $u,v\in\mathcal{C}^{s}$, it holds that
\begin{align}\label{eq:fR-negLip}
        \norm{f_{R}(u)-f_{R}(v)}_{\calC^{-s-\alpha+\epsilon}}\leq C_{R}(1+\norm{u}_{\calC^{s}}^{\lceil s\rceil }+\norm{v}_{\calC^{s}}^{\lceil s\rceil })\norm{u-v}_{\calC^{-s+\epsilon}}.
\end{align}
\end{lemma}
\begin{proof}
Since $G\in C^{\lceil \frac{\sigma-d}{2}\rceil +2}(\R)$ by \cref{ass:G}, we have that $G_{R}\in C^{\lceil \frac{\sigma-d}{2}\rceil +2}_{b}(\R)$ and thus for a constant $C_R>0$, $\norm{\partial^{i} G_{R}}_{\infty}\leq C_R$ for $i=0,\dots,\lceil \frac{\sigma-d}{2}\rceil +2$. Since $f_{R}(u)=\mathcal{D}_{\alpha}G_{R}(u)$, the bound \eqref{eq:fR-bdd} follows easily from \eqref{eq:derivative} and boundedness of $G_R$. Moreover we obtain that by \eqref{eq:derivative} and boundedness of $\partial G_R$,
\begin{align*}
    \norm{f_R(u)-f_{R}(v)}_{\calC^{-\alpha}}&\lesssim \norm{G_R(u)-G_{R}(v)}_{\calC^{0}}\\&\lesssim \norm{G_R(u)-G_{R}(v)}_{L^{\infty}}\\&\leq \norm[\bigg]{\int_{0}^{1} (\partial G_R)(\lambda u+(1-\lambda) v)d\lambda}_{L^{\infty}}\norm{u-v}_{L^{\infty}}\leq C_R \norm{u-v}_{L^{\infty}}
\end{align*}
and thus \eqref{eq:fR-Lip} follows.
The bounds \eqref{eq:fR-lLip} and \eqref{eq:fR-growth} follow in a similar fashion as for \eqref{eq:fR-Lip} replacing $L^{\infty}$ by $\calC^{s}$ and as well as using the following bound: for $i=0,1$, $s\in (0, \lceil\frac{\sigma-d}{2}\rceil)\setminus\N$ and $u\in\calC^{s}$,
\begin{align}\label{eq:G_R-growth}
    \norm{(\partial^{i} G_{R})(u)}_{\calC^{s}}\leq C_{R}(1+\norm{u}_{\calC^{s}}^{\lceil s\rceil}).
\end{align} 
We are left to show \eqref{eq:G_R-growth}. To that aim, let $m=\lceil \frac{\sigma-d}{2}\rceil$. Then \eqref{eq:G_R-growth} can be seen by distinguishing the cases $s\in (0,1), s\in (1,2),\dots, s\in (m-1,m)$.
Let $s\in (r-1,r)$ for any $1\leq r\leq m$. Denote by $\beta=(\beta_1,\dots,\beta_d)\in\N^d_0$ a multi-index and $\partial^{\beta}=\partial^{\beta_1}_1\cdots\partial^{\beta_d}_d$. Then the chain rule (see for example \cite[Lemma 2.3]{Bahouri2011}) gives for a multi-index $\beta$ with $\abs{\beta}=k$, 
\begin{align*}
    \partial^{\beta}[\partial^i G_R(u)]=\sum_{l=1}^{k}\partial^{i+l} G_R(u) \sum_{(\beta^1,\dots,\beta^l)\in\N^l\,:\;\beta^1+\dots+\beta^l=\beta} C_{\beta^1,\dots, \beta^l} \Pi_{j=1}^{l} \partial^{\beta^j} u
\end{align*}
for some nonnegative integers $C_{\beta^1,\dots,\beta^{l}}$. Hence we find for $x,y\in\T^d$,
\begin{align*}
\MoveEqLeft
    \partial^{\beta}[\partial^i G_R(u)](x) - \partial^{\beta}[\partial^i G_R(u)] (y)
    \\&= \sum_{l=1}^{k} \sum_{\beta^1+\dots+\beta^l=\beta} \paren[\Big]{\partial^{i+l} G_R(u(x))-\partial^{i+l} G_R(u(y))}C_{\beta^1,\dots, \beta^l} \Pi_{j=1}^{l} \partial^{\beta^j} u (x) 
    \\&\qquad + \sum_{l=1}^{k} \sum_{\beta^1+\dots+\beta^l=\beta} \partial^{i+l} G_R(u(y))C_{\beta^1,\dots, \beta^l} \paren[\Big]{\Pi_{j=1}^{l} \partial^{\beta^j} u (x) - \Pi_{j=1}^{l} \partial^{\beta^j} u (y)}
\end{align*}
with 
\begin{align*}
    \partial^{i+l} G_R(u(x))-\partial^{i+l} G_R(u(y))= \paren[\Big]{\int_{0}^{1} \partial^{i+l+1} G_R(\lambda u(x)+(1-\lambda)u(y))d\lambda}(u(x)-u(y))
\end{align*}
and
\begin{align*}
\Pi_{j=1}^{l} \partial^{\beta^{j}} u (x) - \Pi_{j=1}^{l} \partial^{\beta^{j}} u (y) = \sum_{j=2}^{l-1}\paren[\Big]{\partial^{\beta^{j}} u (x)- \partial^{\beta^{j}} u (y)}\Pi_{\tilde{j}=1,\dots, j-1} \partial^{\beta^{\tilde{j}}} u (y)\Pi_{\hat{j}=j+1,\dots, l} \partial^{\beta^{\hat{j}}} u (x).
\end{align*}
Thus we have that for $s\in(r-1,r)$ by the equivalence of the Besov and Hölder norm,
\begin{align*}
    \MoveEqLeft
    \norm{(\partial^{i} G_{R})(u)}_{\calC^{s}}\\&\lesssim  \sum_{0\leq\abs{\beta}\leq r-1}\norm{\partial^{\beta} [(\partial^{i} G_{R})(u)]}_{L^\infty}+ \sum_{\abs{\beta}=r-1} \sup_{x\neq y}\frac{\abs{\partial^{\beta}[(\partial^{i} G_{R})(u)](x)-\partial^{\beta}[(\partial^{i} G_{R})(u)](y)}}{\abs{x-y}^{s-r+1}}
    \\&\lesssim C_R (1+\sum_{\abs{\beta}\leq r-1}\norm{\partial^{\beta}u}_{L^\infty}) + C_R \sum_{\abs{\tilde{\beta}}, \abs{\beta'}\leq r-1} \norm{\partial^{\tilde{\beta}}u}_{\calC^{s-r+1}} (1+\norm{\partial^{\beta'}u}_{L^{\infty}}^{r-1})
    \\&\lesssim C_R (1+ \norm{u}_{\calC^{s}}^{r}),
\end{align*}
where we have used that $\norm{\partial^{i+l}G_{R}}_{\infty}\leq C_R$ for $i=0,1$ and any $0\leq l\leq r$, which verifies \eqref{eq:G_R-growth}.

The bound \eqref{eq:fR-negLip} follows from the product estimate \eqref{eq:product} together with \eqref{eq:G_R-growth}, yielding
\begin{align*}
    \norm{f_R(u)-f_{R}(v)}_{\calC^{-s-\alpha+\epsilon}}&\lesssim\norm[\bigg]{\int_{0}^{1} (\partial G_R)(\lambda u+(1-\lambda) v)d\lambda}_{\calC^{s}}\norm{u-v}_{\calC^{-s+\epsilon}} \\&\leq C_R (1+\norm{u}_{\calC^{s}}^{\lceil s\rceil}+\norm{v}_{\calC^{s}}^{\lceil s\rceil})\norm{u-v}_{\calC^{-s+\epsilon}}.
\end{align*}
\end{proof}
\begin{remark}
    From the above proof we extract the following bound on the composition of a function $H\in C^m_b(\R)$ with $u\in \calC^{s}(\T^d)$ with $s\in (0,\infty)\setminus\N$ and $m=\lceil s\rceil$:
    \begin{align}\label{eq:composition}
        \norm{H(u)}_{\calC^{s}}\lesssim \norm{H}_{C^m_b} (1+\norm{u}_{\calC^s}^m).
    \end{align}
\end{remark}

Following the approach introduced in our earlier paper
\cite{DjGK}, we use the stochastic sewing lemma for improving the temporal rate. The stochastic sewing lemma originated from \cite{Khoa}, but we state here a weighted version (see \cite{ABLM, DGL}). The final conclusion of the lemma follows from \cite[Lemma~2.3]{BFG}. Let $$[S,T]_\leq:=\{(s,t)\mid S\leq s<t\leq T\}$$ and $$[S,T]^{\ast}_\leq:=\{(s,t)\mid S\leq s<t< T,\,t-s\leq T-t\}.$$ For a function $\mathcal{A}$ of one variable and $s\leq t$, we write $\mathcal{A}_{st}=\mathcal{A}_{t}-\mathcal{A}_{s}$ and for functions $A$ of two variables and $s\leq u\leq t$, we denote $\delta A_{sut}=A_{st}-A_{su}-A_{ut}$. 
Further, denote by $\E_{s}$ the conditional expectation with respect to $\F_s$.

\begin{lemma}\label{lem:ss}
  Fix $p\geq 2$ and $0\leq S<T\leq 1$. Let $A:[S,T]_{\leq}\to L^{p}(\Omega)$ be such that $A_{st}$ is $\F_{t}$-measurable for all $(s,t)\in[S,T]_{\leq}$. Suppose that there exist $\epsilon_1,\epsilon_2>0,\delta_1,\delta_2\geq 0$ and $C_{1},C_2<\infty$ satisfying $1/2+\epsilon_1-\delta_1>0,\, 1+\epsilon_2-\delta_2>0$ and such that 
  for all $(s,t)\in[S,T]^{\ast}_{\leq}$, $u\in[s,t]$ the bounds
  \begin{align}
      \norm{A_{st}}_{L^{p}(\Omega)}&\leq C_1 (T-t)^{-\delta_1}\abs{t-s}^{1/2+\epsilon_1},\label{eq:SSL-cond1}
      \\
      \norm{\E_{s}[\delta A_{sut}]}_{L^{p}(\Omega)}&\leq C_2 (T-t)^{-\delta_2}\abs{t-s}^{1+\epsilon_2}\label{eq:SSL-cond2}
  \end{align}
  hold. Then there exists a unique $(\F_t)_{t\in[S,T]}$-adapted process $\mathcal{A}:[S,T]\to L^{p}(\Omega)$ such that $\mathcal{A}_{S}=0$ and 
  that there exist $K_1,K_2$ such that for all $(s,t)\in[S,T]^{\ast}_{\leq}$ one has 
  \begin{align}
      \norm{\mathcal{A}_{st}-A_{st}}_{L^{p}(\Omega)}&\leq K_1 (T-t)^{-\delta_1}\abs{t-s}^{1/2+\epsilon_1}+ K_2 (T-t)^{-\delta_2}\abs{t-s}^{1+\epsilon_2},\label{eq:ss1}\\
      \norm{\E_{s}[\mathcal{A}_{st}- A_{st}]}_{L^{p}(\Omega)}&\leq K_2 (T-t)^{-\delta_2}\abs{t-s}^{1+\epsilon_2}\label{eq:ss2}.
  \end{align}
  Furthermore, there exists a constant $C$ depending only on $p,\eps_1,\eps_2$, such that the above bounds hold with $K_{1}=CC_1, K_{2}=C C_2$.
  Finally, there exists a constant $C'$ depending only on $p,\eps_1,\eps_2,\delta_1,\delta_2$, such that for all $(s,t)\in[S,T]_\leq$ one has
  \begin{align}
      \norm{\mathcal{A}_{st}}_{L^{p}(\Omega)}&\leq C'(C_1 \abs{t-s}^{1/2+\epsilon_1-\delta_1}+ C_2 \abs{t-s}^{1+\epsilon_2-\delta_2}).
  \end{align}
\end{lemma}

Let $u^{R}$ be the solution of \eqref{eq:SPDE} with nonlinearity $f_{R}$ from \eqref{eq:trunc-f} in place of $f$, that is
\begin{align}\label{eq:truncated-u}
    u^{R}_t=P_{t}^{\sigma}u_{0}+\int_{0}^{t}P_{t-s}^{\sigma}f_R(u^R_{s})ds+O_t^{\sigma},\quad t\in [0,T].
\end{align}
Let $V^{M,R}$ denote the corresponding accelerated exponential Euler scheme. That is, let again $t_k=\frac{Tk}{M}=k h$ for $k=0,\dots M$ and $V^{M,R}_{0}=u_{0}$ and let
\begin{align*}
    V^{M,R}_{t_{k+1}}=P_{h}^{\sigma}V^{M,R}_{t_{k}} - (-\Delta)^{-\sigma/2}(P_h^{\sigma}-\operatorname{Id})f_R(V^{M,R}_{t_k}) + O_{t_{k+1}}^{\sigma}-P_{h}^{\sigma}O_{t_k}^{\sigma}, \quad k=0,\dots, M-1
\end{align*}
with the convention that $- (-\Delta)^{-\sigma/2}(P_h^{\sigma}-\operatorname{Id})e_0 :=h e_0$ for the zeroth Fourier mode $e_0$.
Written in the mild form, we have that
\begin{align}\label{eq:scheme}
    V^{M,R}_{t}=P_{t}^{\sigma}u_{0}+\int_{0}^{t}P_{t-s}^{\sigma}f_R(V^{M,R}_{k_{M}(s)})ds+O_t^{\sigma},\quad t\in [0,T],
\end{align}
where $k_{M}(s):=\lfloor h^{-1} s\rfloor h$ is the last gridpoint before (or equal to) $s$. The mild form agrees with the inductive form on the time gridpoints $t=t_k$ and is meaningful for any $t\in[0,T]$.\\
We use the following auxiliary lemma to prove a priori estimates for $V^{M,R}$.
\begin{lemma}\label{lem:int-bound}
    Let $v:\{(t_k)\mid k=0,\dots, M\}\to\R^{\geq 0}$ be a function defined on the time grid points $t_k=\frac{Tk}{M}$ of the interval $[0,T]$ for $M\in\N$. Let again $k_{M}(t):=\lfloor h^{-1} t\rfloor h$, $t\in[0,T]$. Let $\gamma\in [0,1)$. Then for any $t\in [0,T]$, it follows that
    \begin{align*}
        \int_{0}^{k_{M}(t)}(k_M(t)-r)^{-\gamma} v_{k_{M}(r)} dr \leq C(\gamma,T) \int_{0}^{t} (t-r)^{-\gamma}v_{k_{M}(r)} dr.
    \end{align*}
\end{lemma}
\begin{proof}
Let for $t\in [0,T]$, $l\in \{0,\dots,M\}$ be such that $k_{M}(t)=t_l$. If $l=0$, the claim is trivial, so that we can assume $l\geq 1$.
We first prove that for $k=0,\dots, l-1$, there exists $C(\gamma,T)>0$, such that
\begin{equ}\label{eq:int-claim}
     \int_{t_{k}}^{t_{k+1}}(t_l -r)^{-\gamma} dr \leq C(\gamma,T)  \int_{t_{k}}^{t_{k+1}}(t -r)^{-\gamma} dr.
\end{equ}
Let first $k=l-1$, in which case  we recall that $t\leq t_l+TM^{-1}$, such that $(t-r)\leq 2 TM^{-1}$ and thus
\begin{equ}
     \int_{t_{l-1}}^{t_{l}}(t_l -r)^{-\gamma} dr=\frac{M^{\gamma-1} T^{1-\gamma} }{1-\gamma}, \qquad \int_{t_{l-1}}^{t_{l}}(t -r)^{-\gamma} dr\geq TM^{-1}(2 TM^{-1})^{-\gamma}=\frac{M^{\gamma-1}T^{1-\gamma}}{2^\gamma},
\end{equ}
which implies that \eqref{eq:int-claim} holds in this case. 
Moreover, if $k<l-1$, we have that $(t_l-r)\geq \frac{1}{3}(t-r) $ for $r \in [t_k, t_{k+1}]$, so that \eqref{eq:int-claim} follows by comparing the integrands. Using \eqref{eq:int-claim}, we then obtain
    \begin{align*}
        \int_{0}^{k_{M}(t)}(k_M(t)-r)^{-\gamma} v_{k_{M}(r)} dr&= \sum_{k=0}^{l-1} v_{t_{k}} \int_{t_{k}}^{t_{k+1}}(t_l -r)^{-\gamma} dr 
        \\&\leq C(\gamma,T) \sum_{k=0}^{l-1} v_{t_{k}} \int_{t_{k}}^{t_{k+1}}(t -r)^{-\gamma} dr 
        \\&\leq C(\gamma,T)  \int_{0}^{t} (t-r)^{-\gamma}  v_{k_{M}(r)} dr  ,
    \end{align*}
    where the constant $C(\gamma,T)$ changes from inequality to inequality.
\end{proof}
Next, we establish the a priori estimate for the numerical solution $V^{M,R}$.

 \begin{lemma}\label{lem:remainder}
 Let $R>0$ and $V^{M,R}$ be the accelerated exponential Euler scheme as in \eqref{eq:scheme}.
 Then the following a priori estimate holds true. For any $\epsilon\in (0,\frac{\sigma-d}{2})$, $p\geq 1$ and $\lambda\in [0,\frac{\sigma}{2})$,
 \begin{align}\label{eq:VRM-moment}
\sup_{M\in\N}\E\norm{V^{M,R}}_{C_T^{\lambda/\sigma}\calC^{\frac{\sigma-d}{2}-\lambda-\epsilon}}^p\leq C(\mathcal{M}(\epsilon),p,\epsilon,\sigma,\lambda, \alpha, T)<\infty. 
 \end{align}
 \end{lemma}
 \begin{proof}
  Let $p\geq 1$ and $t_k=\frac{Tk}{M}$ for $k=0,\dots, M$. We first prove by induction that $\E[\norm{V^{M,R}_{t_k}}_{L^{\infty}}^p<\infty$ for all $k=0,\dots, M$. Since $V^{M,R}_{t_0}=u_0\in L^{\infty}$, we may use \eqref{eq:fR-bdd} and $\sigma>\alpha$ together with $\calC^{\eps}\hookrightarrow L^{\infty}$, to infer that for small $\eps>0$,
   \begin{align*}
     \norm{V^{M,R}_{t_1}}_{L^{\infty}}&\lesssim \norm{u_{0}}_{L^{\infty}}+\int_{0}^{t_1}(t_1-r)^{-\frac{\alpha+\eps}{\sigma}}\norm{f_R(V^{M,R}_{t_0})}_{\calC^{-\alpha}}dr + \norm{O^{\sigma}}_{C_T L^{\infty}}
     \\&\lesssim \norm{u_{0}}_{L^{\infty}}+C_R\int_{0}^{t_1}(t_1-r)^{-\frac{\alpha+\eps}{\sigma}}dr 
     + \norm{O^{\sigma}}_{C_T L^{\infty}}
      \\&\lesssim \norm{u_{0}}_{L^{\infty}}+C_R T^{1-\frac{\alpha+\eps}{\sigma}}
      + \norm{O^{\sigma}}_{C_T L^{\infty}}.
   \end{align*}
   This implies that $\E\norm{V^{M,R}_{t_1}}_{L^{\infty}}^p<\infty$ and in particular that almost surely $V^{M,R}_{t_1}\in L^{\infty}$. Now let $k\in\{2,\dots, M\}$ and assume that $\E\norm{V^{M,R}_{t_l}}_{L^{\infty}}^p<\infty$ for all $l\leq k-1$. Then we can show that $\E\norm{V^{M,R}_{t_{k}}}_{L^{\infty}}^p<\infty$ by the same argument as for the step $t_0\to t_1$ above using again \eqref{eq:fR-bdd} together with almost surely, for all $l\leq k-1$, $V^{M,R}_{t_l}\in L^{\infty}$.  
   Furthermore by a similar bound for $V^{M,R}_{t_k}$ as for $V^{M,R}_{t_1}$ above, we obtain that
   \begin{align*}
       \E[\max_k\norm{V^{M,R}_{t_k}}_{L^{\infty}}^p]\lesssim \E[ (\norm{u_{0}}_{L^{\infty}}+\norm{O^{\sigma}}_{C_T L^{\infty}}+C_R T^{1-\frac{\alpha+\eps}{\sigma}})^p]<\infty.
   \end{align*}
Next we infer the time-regularity bound. We have that for $0\leq s\leq t\leq T$, using again \eqref{eq:fR-bdd} with $V^{M,R}_{t_k}\in L^{\infty}$ for all $k=0,\dots,M$ almost surely, and the semigroup estimates \eqref{eq:HK-1}, as well as $\sigma> 2\alpha-d$, $\lambda\in [0,\frac{\sigma}{2})$, $\epsilon\in (0,\frac{\sigma-d}{2})$, almost surely
    \begin{align}\label{eq:VRM-reg}
    \MoveEqLeft
     \norm{V^{M,R}_{t}-V^{M,R}_{s}}_{\calC^{\frac{\sigma-d}{2}-\lambda-\epsilon}}\nonumber\\&\lesssim \norm{u_{0}}_{\calC^{\frac{\sigma-d}{2}-\epsilon}}+\int_{s}^{t}(t-r)^{\min(-\frac{2\alpha+\sigma-d}{2\sigma}+\frac{\epsilon}{\sigma}+\frac{\lambda}{\sigma},0)}\norm{f_R(V^{M,R}_{k_{M}(r)})}_{\calC^{-\alpha}}dr \nonumber\\&\qquad+ \int_{0}^{s} (t-s)^{\frac{\lambda}{\sigma}}(s-r)^{-\frac{2\alpha+\sigma-d}{2\sigma}+\frac{\epsilon}{\sigma}}\norm{f_R(V^{M,R}_{k_{M}(r)})}_{\calC^{-\alpha}}dr\nonumber
     \\&\qquad + \norm{O^{\sigma}}_{C_T^{\frac{\lambda}{\sigma}}\calC^{\frac{\sigma-d}{2}-\epsilon-\lambda}} (t-s)^{\frac{\lambda}{\sigma}}
     \nonumber\\&\lesssim \norm{u_{0}}_{\calC^{\frac{\sigma-d}{2}-\epsilon}}+C_R  (t-s)^{\frac{\lambda}{\sigma}}
     +  \norm{O^{\sigma}}_{C_T^{\frac{\lambda}{\sigma}}\calC^{\frac{\sigma-d}{2}-\epsilon-\lambda}} (t-s)^{\frac{\lambda}{\sigma}}
   \end{align}
   and thus the claim follows from \cref{prop:ou}.
 \end{proof}

We now state the following auxiliary results in preparation for the proof of the main theorem.
\begin{lemma}\label{lem:conv-R}
Let for $R>0$, $u^R$ and $V^{M,R}$ be defined as in \eqref{eq:truncated-u} and \eqref{eq:scheme}, respectively. Then $u_R$ is the unique mild solution to \eqref{eq:truncated-u}.
Furthermore, the following moment bound holds true for any $\epsilon\in (0,\frac{\sigma-d}{2})$, $\lambda\in [0,\frac{\sigma}{2})$ and $p\geq 1$,
\begin{align}\label{eq:uR-moment}
    \E[\norm{u^R}_{C_{T}^{\frac{\lambda}{\sigma}}\calC^{\frac{\sigma-d}{2}-\lambda-\epsilon}}^p]<\infty
\end{align}
and for $d^{R}_t:=u^{R}_t-P_{t}^{\sigma}u_0-O^{\sigma}_t$, $t\in[0,T]$,
we have that, for any $\epsilon\in (0,\frac{\sigma-d}{2}), \lambda\in [0,\sigma)$ and $p\geq 1$,
\begin{align}\label{eq:dR-moment}
    \E[\norm{d^R}_{C_{T}^{\frac{\lambda}{\sigma}}\calC^{\frac{3\sigma-d}{2}-\alpha-\lambda-\epsilon}}^p]<\infty.
\end{align}
 Let $$\tau_{R}:=\inf\{t\in[0,T]\mid \norm{u^{R}_t}_{L^{\infty}}\geq R\}$$ with $\inf\emptyset=+\infty$. Then we have that $u^{R}_t=u_t$ for all $t\leq \tau_{R}\wedge T$. Let $e^{M,R}:=u^R-V^{M,R}$ be the error. Then the following convergence hold true for any $p\geq 1$,
\begin{align}\label{eq:L2-conv-R}
    \lim_{M\to\infty}\E\sup_{t\in[0,T]}\norm{e^{M,R}_t}_{L^{\infty}}^p=0.
\end{align}
\end{lemma}
\begin{proof}
The wellposedness of the equation for $u^R$ is immediate, but we give the argument, since $f_R$ decreases regularity. To that aim, for $p\geq 1$, we define the fixed point map $$\Phi: L^p(\Omega, C_T L^{\infty})\ni u\mapsto \Phi(u)=\paren[\Big]{P_{t}^{\sigma} u_0+\int_{0}^{t}P^{\sigma}_{t-r}f_R(u_r)dr+O^{\sigma}_t}_{t\in[0,T]}\in L^p(\Omega, C_T L^{\infty})$$ and first check that $\Phi$ indeed maps into $L^p(\Omega, C_TL^{\infty})$.
Since $\sigma>\alpha$ and using \eqref{eq:fR-bdd} with $u\in C_T L^{\infty}$ almost surely, we have that by $\calC^{\epsilon}\hookrightarrow L^{\infty}$ and the semigroup estimates \eqref{eq:HK-1}, almost surely for all $t\in[0,T]$,
   \begin{align}\label{eq:b}
     \norm{\Phi(u)_t}_{L^{\infty}}&\lesssim \norm{u_{0}}_{L^{\infty}}+\int_{0}^{t}(t-r)^{-\frac{\alpha+\eps}{\sigma}}\norm{f_R(u_r)}_{\calC^{-\alpha}}dr + \norm{O^{\sigma}}_{C_T L^{\infty}}\nonumber
     \\&\lesssim \norm{u_{0}}_{L^{\infty}}+C_R\int_{0}^{t}(t-r)^{-\frac{\alpha+\eps}{\sigma}}dr 
     + \norm{O^{\sigma}}_{C_T L^{\infty}}.
   \end{align} 
Thus we obtain that for small enough $\epsilon>0$,
   \begin{align}\label{eq:lp}
       \E[\sup_{t\in[0,T]}\norm{\Phi(u)_t}_{L^{\infty}}^p]\lesssim \E[ (\norm{u_{0}}_{L^{\infty}}+\norm{O^{\sigma}}_{C_T L^{\infty}}+C_R T^{1-\frac{\alpha+\eps}{\sigma}})^p]<\infty.
   \end{align}
Continuity in time follows from the continuity of the time integral and $O\in L^p(\Omega,C_T L^{\infty})$.
Next, we argue that for small enough $T$, $\Phi$ is a contraction and thus a unique fixed point exists. But this is clear by an analogue estimate as above using the global Lipschitz bound of $f_R$ into $\calC^{-\alpha}$, that is \eqref{eq:fR-Lip} in place of \eqref{eq:fR-bdd}.
Patching the solutions on subintervals together, we thus arrive at a unique mild solution $u^{R}$ satisfying the bound \eqref{eq:lp}.
    Furthermore, since $f(u)=f_{R}(u)$ for all continuous $u$ with $\norm{u}_{L^{\infty}}=\norm{u}_{\infty}\leq R$ by definition \eqref{eq:trunc-f} and we have that by \cref{ass:1} $x\mapsto u_t(x)$ is continuous, we obtain that $f_R(u_t)=f(u_t)$ for $t\leq\tau_R\wedge  T$ and thus that $(u_t)_{t\in[0,\tau_R\wedge T]}$ solves the equation for $u^R$, such that by uniqueness $u^R_t=u_t$ for $t\leq \tau_R\wedge T$.\\
    Next, we prove \eqref{eq:uR-moment} and \eqref{eq:dR-moment}. 
    Due to the growth bound \eqref{eq:fR-bdd} on $f_R$ and the semigroup estimates \eqref{eq:HK-1}, we see that almost surely for any $t\in[0,T]$,
    \begin{align}\label{eq:lambda0}
     \norm{u^{R}_t}_{\calC^{\frac{\sigma-d}{2}-\epsilon}}&\lesssim \norm{u_{0}}_{\calC^{\frac{\sigma-d}{2}-\epsilon}}+\int_{0}^{t}(t-r)^{\frac{d-\sigma-2\alpha+2\epsilon}{2\sigma}}\norm{f_R(u^{R}_t)}_{\calC^{-\alpha}}dr + \norm{O^{\sigma}}_{C_T \calC^{\frac{\sigma-d}{2}-\epsilon}}\nonumber
     \\&\lesssim \norm{u_{0}}_{\calC^{\frac{\sigma-d}{2}-\epsilon}}+C_R T^{\frac{d+\sigma-2\alpha+2\epsilon}{2\sigma}}
     + \norm{O^{\sigma}}_{C_T \calC^{\frac{\sigma-d}{2}-\epsilon}}.
    \end{align}
    Thus the bounds on $O^{\sigma}$ from \cref{prop:ou} imply that \eqref{eq:uR-moment} holds for $\lambda=0$. 
Next, we note that $u_{t}=P_{t}^{\sigma} u_0 + d^{R}_{t} + O_{t}^{\sigma}$ and that the bound \eqref{eq:uR-moment} holds for $u^{R}$ replaced by $O^{R}$ due to \cref{prop:ou}. Moreover by the semigroup estimates for $0\leq s\leq t\leq T$ and $\lambda\in [0,\sigma]$,
    \begin{align}\label{eq:bound-u0}
        \norm{P_{t}^{\sigma} u_0 - P_{s}^{\sigma}u_0}_{\calC^{\frac{\sigma-d}{2}-\epsilon-\lambda}} \lesssim \abs{t-s}^{\lambda/\sigma} \norm{u_0}_{\calC^{\frac{\sigma-d}{2}-\epsilon}}.
    \end{align}
    Thus, we are left to prove the bound \eqref{eq:dR-moment} on $d^R$, which in particular implies the bound \eqref{eq:uR-moment} on $u^R$ together with \eqref{eq:bound-u0} and the representation above, since $\sigma>\alpha$.
     We have that for $\lambda\in [0,\sigma)$ by the semigroup estimates \eqref{eq:HK-1} and \eqref{eq:fR-growth},  for $0\leq s\leq t\leq T$,
 \begin{align}\label{eq:d-reg}
 \MoveEqLeft
     \norm{d_t^R - d_s^R}_{\calC^{\frac{3\sigma-d}{2}-\alpha-\lambda-\epsilon}}\nonumber
     \\&\leq \int_{s}^{t}\norm{P_{t-r}^{\sigma} f_R(u_r^R)}_{\calC^{\frac{3\sigma-d}{2}-\alpha-\lambda-\epsilon}}dr + \int_0^s \norm{P_{s-r}^{\sigma}(\operatorname{Id}-P_{t-s}^{\sigma})f_R(u_r^R)}_{\calC^{\frac{3\sigma-d}{2}-\alpha-\lambda-\epsilon}}dr\nonumber
     \\&\lesssim \int_{s}^{t} (t-r)^{-1+\frac{\lambda}{\sigma}}\norm{f_R(u_r^R)}_{\calC^{\frac{\sigma-d}{2}-\alpha-\epsilon}}dr \nonumber
     \\&\qquad+ \int_0^s (t-s)^{\frac{\lambda}{\sigma}}(s-r)^{-1+\frac{\epsilon}{2\sigma}}\norm{f_R(u_r^R)}_{\calC^{\frac{\sigma-d}{2}-\alpha-\frac{\epsilon}{2}}} dr \nonumber
     \\&\lesssim (t-s)^{\frac{\lambda}{\sigma}}C_R(1+\norm{u^R}_{C_T\calC^{\frac{\sigma-d-\epsilon}{2}}}^m)
 \end{align}
 where $m=\lceil \frac{\sigma-d-\epsilon}{2}\rceil$ 
 and thus \eqref{eq:dR-moment} follows from taking the $L^p$-moment of \eqref{eq:lambda0}.
   It remains to prove \eqref{eq:L2-conv-R}.
    Using the semigroup estimates with $\sigma>\alpha$ and the Lipschitz bound on $f_R$ from \eqref{eq:fR-Lip}, as well as $\calC^{\epsilon}\hookrightarrow L^{\infty}$, we have that for $\epsilon\in (0,\sigma-\alpha)$
     \begin{align}\label{eq:error-b}
        \norm{V^{M,R}_t-u^R_t}_{L^{\infty}}&= \norm[\bigg]{\int_{0}^{t}P_{t-s}^{\sigma}[f_R(V^{M,R}_{k_{M}(r)})-f_{R}(u^R_r)]dr}_{L^{\infty}}\nonumber
        \\&\leq \norm[\bigg]{\int_{0}^{t}P_{t-r}^{\sigma}[f_R(V^{M,R}_{k_{M}(r)})-f_{R}(u^R_{k_{M}(r)}]dr}_{L^{\infty}} \nonumber
        \\&\quad + \norm[\bigg]{\int_{0}^{t}P_{t-r}^{\sigma}[f_R(u^R_{k_{M}(r)})-f_{R}(u^R_{r})]dr}_{L^{\infty}}\nonumber
        \\& \lesssim C_R \int_{0}^{t}(t-r)^{-\frac{\alpha+\epsilon}{\sigma}} \norm{V_{k_M(r)}^{M,R}-u^R_{k_M (r)}}_{L^{\infty}}dr \nonumber
        \\&\quad + C_R \int_{0}^{t}(t-r)^{-\frac{\alpha+\epsilon}{\sigma}} \norm{u_{k_M(r)}^{R}-u^R_{r}}_{L^{\infty}}dr \nonumber
        \\& \lesssim C_R \int_{0}^{t}(t-r)^{-\frac{\alpha+\epsilon}{\sigma}} \norm{V_{k_M(r)}^{M,R}-u^R_{k_M (r)}}_{L^{\infty}}dr \nonumber
         \\&\quad + C_R T^{\frac{\sigma-\alpha-\epsilon}{\sigma}}M^{-\frac{\sigma-d}{2\sigma}+\epsilon}\norm{u^R}_{C_T^{\frac{\sigma-d}{2\sigma}-\epsilon\sigma}L^{\infty}}.
     \end{align}  
     Plugging $k_{M}(t)$ in place of $t$ in \eqref{eq:error-b} and using \cref{lem:int-bound}, we see that
     \begin{align*}
        \norm{V^{M,R}_{k_{M}(t)}-u^R_{k_{M}(t)}}_{L^{\infty}}&\leq  C_R \int_{0}^{k_{M}(t)}(k_{M}(t)-r)^{-\frac{\alpha+\epsilon}{\sigma}} \norm{V_{k_M(r)}^{M,R}-u^R_{k_M (r)}}_{L^{\infty}}dr \nonumber
         \\&\quad + C_R T^{\frac{\sigma-\alpha-\epsilon}{\sigma}}M^{-\frac{\sigma-d}{2\sigma}+\epsilon}\norm{u^R}_{C_T^{\frac{\sigma-d}{2\sigma}-\epsilon}L^{\infty}}
         \\&\leq  C_R \int_{0}^{t}(t-r)^{-\frac{\alpha+\epsilon}{\sigma}} \norm{V_{k_M(r)}^{M,R}-u^R_{k_M (r)}}_{L^{\infty}}dr \nonumber
         \\&\quad + C_R T^{\frac{\sigma-\alpha-\epsilon}{\sigma}}M^{-\frac{\sigma-d}{2\sigma}+\epsilon}\norm{u^R}_{C_T^{\frac{\sigma-d}{2\sigma}-\epsilon}L^{\infty}}.
     \end{align*}
     Applying the Henry-Gr\"onwall inequality (\cite[Lemma 7.1.1]{Henry}) for $y(t)=\norm{V^{M,R}_{k_{M}(t)}-u^R_{k_{M}(t)}}_{L^{\infty}}$, $0\leq t\leq T$, and taking $p$-th moments, we obtain that
     \begin{align*}
         (\E\sup_{k=0,\dots, M} \norm{V^{M,R}_{t_k}-u^R_{t_k}}_{L^{\infty}}^p)^{1/p}\lesssim C_R T^{\frac{\sigma-\alpha}{\sigma}}M^{-\frac{\sigma-d}{2\sigma}+\epsilon}(\E\norm{u^R}_{C_T^{\frac{\sigma-d}{2\sigma}-\epsilon}L^{\infty}}^p)^{1/p}.
     \end{align*}
     From the time regularities of $V^{M,R}$ and $u^{R}$ from \eqref{eq:VRM-moment} and \eqref{eq:uR-moment} and triangle inequality, we see that for $t\in[0,T]$, $$\norm{V^{M,R}_{t}-u^R_{t}}_{L^{\infty}}\leq \norm{V^{M,R}_{k_{M}(t)}-u^R_{k_{M}(t)}}_{L^{\infty}} + M^{-\frac{\sigma-d}{2\sigma}+\epsilon}[\norm{V^{M,R}}_{C_{T}^{\frac{\sigma-d}{2\sigma}-\epsilon}L^{\infty}}+ \norm{u^R}_{C_{T}^{\frac{\sigma-d}{2\sigma}-\epsilon}L^{\infty}}]$$ 
     and thus together with the above, it follows that as $M\to\infty$, $$\E[\sup_{t\in[0,T]}\norm{e^{M,R}_t}_{L^{\infty}}^p]\to 0.$$
     We note that in \cref{prop:main} below we even prove rate $M^{-\frac{3}{2}+\frac{d+2\alpha}{2\sigma}+\epsilon}$ for the latter convergence. Here the convergence with some positive rate is enough to deduce \eqref{eq:L2-conv-R}.\\
\end{proof}

The following lemma is a slight adaptation of \cite[Lemma 4.8]{Printems}.
\begin{lemma}\label{lem:conv}
Let $u$ and $V^{M}$ be defined as in \eqref{eq:mild} and \eqref{eq:scheme1}, respectively. Let $e^{M}=u-V^{M}$ be the error. Then we have that for any $\epsilon\in (0,1)$
\begin{align}\label{eq:pconv}
   \limsup_{M\to\infty}\p(\sup_{t\in[0,T]}\norm{e^{M}_t}_{L^{\infty}}> \epsilon)=0,
\end{align}
i.e. the scheme converges in probability. Moreover we have that
\begin{align}\label{eq:a-priori-num}
\lim_{R\to\infty}\limsup_{M\to\infty}\p(\sup_{t\in[0,T]}\norm{V^M_t}_{L^{\infty}}\geq R)=0.
\end{align}
\end{lemma}
\begin{proof}
The proof of the lemma follows in a similar fashion as the proof of \cite[Lemma 4.8]{Printems}. We repeat the arguments in order to adjust the proof to our setting.  By \cref{ass:1} we find that in particular,
 \begin{align}\label{eq:a-priori}
\lim_{L\to\infty}\p(\sup_{t\in[0,T]}\norm{u_t}_{L^{\infty}} \geq L) = 0.
 \end{align}
Let for $R>1$, 
$$\tau_{R-1}:=\inf\{t\in[0,T]\mid \norm{u_t}_{L^{\infty}}\geq R-1\}$$
and let for $\epsilon\in (0,1)$,
$$t_{\epsilon}:=\inf\{t\in (0,T]\mid \norm{e^{M}_{t}}_{L^{\infty}}\geq \epsilon\}$$
with the infimum over an empty set to be defined as $+\infty$.
Then we have that
\begin{align*}
    \p(\sup_{t\in[0,T]}\norm{e^{M}_{t}}>\epsilon)&\leq \p( \sup_{t\in[0,T]}\norm{e^{M}_{t}}_{L^{\infty}}>\epsilon, \tau_{R-1}\leq T) + \p(\sup_{t\in[0,T]}\norm{e^{M}_{t}}_{L^{\infty}}>\epsilon, \tau_{R-1}>T, t_{\epsilon}\leq T) 
    \\&\quad + \p(\sup_{t\in[0,T]}\norm{e^{M}_{t}}_{L^{\infty}}>\epsilon, \tau_{R-1}>T, t_{\epsilon}>T)
    \\&\leq \p(\sup_{t\in[0,T]}\norm{u_t}_{\infty}\geq R-1) + \p(\sup_{t\in[0,T]}\norm{e^{M}_{t}}_{L^{\infty}}>\epsilon, \tau_{R-1}>T, t_{\epsilon}\leq T),
\end{align*}
where we have used that the event $\{\sup_{t\in[0,T]}\norm{e^{M}_{t}}_{L^{\infty}}>\epsilon, \tau_{R-1}>T, t_{\epsilon}>T\}$ is empty.
 On the event $\{\tau_{R-1}> T,\, t_{\epsilon}\leq T\}$, we have that by \cref{lem:conv-R}, $u^{R}_t=u_t$ 
 for all $t\in[0,T]$ since $\tau_{R-1}> T$. Furthermore on that event, we infer that $\norm{V^{M}_{r}}_{L^{\infty}}\leq \norm{V^{M}_{r}-u_{r}}_{L^{\infty}}+\norm{u_{r}}_{L^{\infty}}\leq \epsilon + R-1\leq R $ for $r\in[0,t_{\epsilon}]$. Thus it follows that $V^{M ,R}_r=V^{M}_r$ for $r\in[0,t_{\epsilon}]$,
 which can be seen using the definition of $V^M$ with nonlinearity $f$ which satisfies $f(v)=f_R(v)$ for all continuous $v$ with $\norm{v}_{\infty}=\norm{v}_{L^{\infty}}\leq R$ and by iterating over the intervals $[0,t_1], \dots, [t_{k-1}, t_{\epsilon}]$, where $k$ is such that $t_{k-1}< t_{\epsilon}\leq t_k$ and $t_{l}=\frac{Tl}{h}$ are the time grid points. Hence, on that event we have that $e^{M ,R}_r=e^{M}_r$ for all $r\in [0,t_\epsilon]$ and thus $\sup_{t\in[0,T]}\norm{e^{M,R}_{t}}_{L^{\infty}}\geq \norm{e^{M,R}_{t_{\epsilon}}}_{L^{\infty}}=\norm{e^{M}_{t_{\epsilon}}}_{L^{\infty}}\geq \epsilon>\epsilon/2$.
 Thus together we deduce that
 \begin{align*}
     \p(\sup_{t\in[0,T]}\norm{e^{M}_{t}}>\epsilon)&\leq \p(\sup_{t\in[0,T]}\norm{u_t}_{L^{\infty}}\geq R-1) + \p(\sup_{t\in[0,T]}\norm{e^{M,R}_{t}}_{L^{\infty}}>\epsilon/2).
 \end{align*}
Hence the claim \eqref{eq:pconv} follows from \eqref{eq:L2-conv-R} together with Markov's inequality and \eqref{eq:a-priori}. To see \eqref{eq:a-priori-num}, we observe that for $R>1$ and $\epsilon\in (0,1)$ by triangle inequality,
 \begin{align*}
     \p(\sup_{t\in[0,T]}\norm{V^{M}_t}_{L^{\infty}}\geq R)\leq \p(\sup_{t\in[0,T]}\norm{u_t}_{L^{\infty}}\geq R-\epsilon) + \p(\sup_{t\in[0,T]}\norm{e^{M}_t}_{L^{\infty}}\geq \epsilon)
 \end{align*}
 and thus the claim follows from \eqref{eq:pconv} and \eqref{eq:a-priori}.
\end{proof}

\section{Proof of the main \cref{thm:main}}\label{sec:proof}
Recall that we assume \cref{ass:G} and \cref{ass:1}.
Let again $e^M =u-V^M$ and $e^{M,R}=u^{R}-V^{M,R}$. Define for $R>0$ the stopping time (with $\inf\emptyset=+\infty$)
\begin{align}
   \tau_{R}^M=\inf\{t\leq T\mid \norm{u_{t}}_{L^{\infty}}\geq R\}\wedge  \inf\{t\leq T\mid \norm{V^M_{t}}_{L^{\infty}}\geq R\}.
\end{align}
Using that $e^{M}=e^{M,R}$ on the event that $\tau_{R}^M>T$, which holds by the uniqueness of solutions to the equation for $u^R$ by \cref{lem:conv-R} and the definition of $V^M$, since $f(u)=f_R(u)$ for all $u$ with $\norm{u}_{\infty}\leq R$, we obtain
\begin{align} 
\MoveEqLeft
    \p(\sup_{t\in[0,T]}\norm{e^M_t}_{L^{\infty}}\geq \delta M^{-\frac{3}{2}+\frac{d+2\alpha}{2\sigma}+\epsilon})\nonumber
    \\&\leq \p(\tau_{R}^M >T ,\sup_{t\in[0,T]}\norm{e^M_t}_{L^{\infty}}\geq \delta M^{-\frac{3}{2}+\frac{d+2\alpha}{2\sigma}+\epsilon})
    + \p(\tau_{R}^M\leq T,\sup_{t\in[0,T]}\norm{e^M_t}_{L^{\infty}}\geq \delta M^{-\frac{3}{2}+\frac{d+2\alpha}{2\sigma}+\epsilon})\nonumber
    \\&\leq \p(\sup_{t\in[0,T]}\norm{e^{M,R}_t}_{L^{\infty}}\geq \delta M^{-\frac{3}{2}+\frac{d+2\alpha}{2\sigma}+\epsilon})+\p(\sup_{t\in[0,T]}\norm{u_t}_{L^{\infty}}\geq R)+\p(\sup_{t\in[0,T]}\norm{V^M_t}_{L^{\infty}}\geq R).\label{eq:decomp}
\end{align}
By \cref{ass:1}, we observe that 
$$\lim_{R\to\infty}\p(\sup_{t\in[0,T]}\norm{u_t}_{L^{\infty}}\geq R)=0.$$ Due to \cref{lem:conv}, we see that $$\lim_{R\to\infty}\limsup_{M\to\infty}\p(\sup_{t\in[0,T]}\norm{V^M_t}_{L^{\infty}}\geq R)=0.$$ In order to prove \cref{thm:main}, it is left to prove that for each fixed $R>0$, 
\begin{align}\label{eq:error-rate}
\limsup_{M\to\infty}\p(\sup_{t\in[0,T]}\norm{e^{M,R}_t}_{L^{\infty}}\geq \delta M^{-\frac{3}{2}+\frac{d+2\alpha}{2\sigma}+\epsilon})=0.
\end{align}
Indeed with \eqref{eq:error-rate} and \eqref{eq:decomp}, we obtain that
\begin{align*}
\MoveEqLeft
    \limsup_{M\to\infty}\p(\sup_{t\in[0,T]}\norm{e^M_t}_{L^{\infty}}\geq \delta M^{-\frac{3}{2}+\frac{d+2\alpha}{2\sigma}+\epsilon})
    \\&\leq\p(\sup_{t\in[0,T]}\norm{u_t}_{L^{\infty}}\geq R)+\limsup_{M\to\infty}\p(\sup_{t\in[0,T]}\norm{V^M_t}_{L^{\infty}}\geq R)
\end{align*} and by sending $R\to\infty$, \cref{thm:main} follows.\\

The convergence \eqref{eq:error-rate} follows from Markov's inequality and \cref{prop:main} below 
for an even smaller choice of $\epsilon>0$. In \cref{prop:main} below we prove a rate of convergence for the scheme $V^{M,R}$ in $L^p(\Omega)$. Thus we observe that the convergence rate in probability of the scheme $V^{M}$ is determined by the convergence rate in $L^{p}(\Omega)$ of the scheme $V^{M,R}$.\\
\cref{prop:main} follows from \cref{prop:sewing} below in which an approach based on the stochastic sewing lemma is used following a similar strategy as in the proof of \cite[Proposition 6.5]{DjGK}.  
The main difference compared to that other proof is that we lose in the rate due to the assumption that $f_R(u)$ admits $\alpha$ units less regularity compared to $u$ itself. Heuristically this means that we have to adapt the proof of \cite[Proposition 6.5]{DjGK} by replacing the $L^{\infty}$-norm in which the error is measured by the stronger $\calC^{\alpha}$-norm. 
The proof of \cref{prop:sewing} is given after the proof of \cref{prop:main} below.
\begin{proposition}\label{prop:sewing}
    Let $T>0$, $p\geq 1$. Let $R\in (0,\infty)$. Then we have that for any $\epsilon\in (0,\frac{\sigma-d}{2})$, there exists $C=C(T,\epsilon, p,R)$ such that
    \begin{align}\label{eq:sewing}
    \E\bigg[\sup_{t\in[0,T]}\norm[\bigg]{\int_{0}^{t}P_{t-s}^{\sigma}[f_{R}(u^{R}_s)-f_{R}(u^{R}_{k_{M}(s)})]ds}_{L^{\infty}}^p\bigg]^{1/p}\leq C M^{-\frac{3}{2}+\frac{d+2\alpha}{2\sigma}+\epsilon}.
\end{align}
\end{proposition}
\begin{proposition}\label{prop:main}
   Let $T>0$, $p\geq 1$. Let $R\in (0,\infty)$.   
Then we have that for any $\epsilon\in (0,\frac{\sigma-d}{2})$, there exists $C=C(T,\epsilon, p,R)$ such that 
    \begin{align}\label{eq:error-R-Linfty}
    \E[\sup_{t\in[0,T]}\norm{e^{M,R}_t}_{L^{\infty}}^p]^{1/p}\leq C M^{-\frac{3}{2}+\frac{d+2\alpha}{2\sigma}+\epsilon}.
    \end{align}
\end{proposition}
\begin{proof}
We have that 
\begin{align*}
 e^{M,R}_t= u^{R}-V^{M,R}=\int_{0}^{t}P_{t-s}^{\sigma}[f_{R}(u^{R}_{s})-f_{R}(u^{R}_{k_{M}(s)})]ds + \int_{0}^{t}P_{t-s}^{\sigma}[f_{R}(u^{R}_{k_{M}(s)})-f_R(V^{M,R}_{k_{M}(s)})]ds.
\end{align*}
Due to the Lipschitz bound \eqref{eq:fR-Lip}, semigroup estimates \eqref{eq:HK-1}, $\calC^{\epsilon}\hookrightarrow L^{\infty}$ and using \cref{lem:int-bound}, we see that for small $\epsilon>0$
\begin{align*}
  \norm{e^{M,R}_{k_{M}(t)}}_{L^{\infty}}&\lesssim \sup_{r\in[0,T]}\norm[\bigg]{\int_{0}^{r}P_{r-s}^{\sigma}[f_{R}(u^{R}_s)-f_{R}(u^{R}_{k_{M}(s)})]ds}_{L^{\infty}}  \\&\qquad+ C_{R}\int_{0}^{k_{M}(t)} (k_{M}(t)-s)^{-\frac{\alpha+\epsilon}{\sigma}} \norm{e^{M,R}_{k_{M}(s)}}_{L^{\infty}} ds
  \\&\lesssim \sup_{r\in[0,T]}\norm[\bigg]{\int_{0}^{r}P_{r-s}^{\sigma}[f_{R}(u^{R}_s)-f_{R}(u^{R}_{k_{M}(s)})]ds}_{L^{\infty}}  \\&\qquad+ C_{R}\int_{0}^{t} (t-s)^{-\frac{\alpha+\epsilon}{\sigma}} \norm{e^{M,R}_{k_{M}(s)}}_{L^{\infty}} ds
\end{align*} and thus by the Henry-Gr\"onwall inequality (\cite[Lemma 7.1.1]{Henry}) applied to $y(t):=\norm{e^{M,R}_{k_{M}(t)}}_{L^{\infty}}$, $t\in[0,T]$, we obtain that 
\begin{align}\label{eq:Gb}
  \max_{k}\norm{e^{M,R}_{t_k}}_{L^{\infty}}=\sup_{t\in[0,T]}\norm{e^{M,R}_{k_{M}(t)}}_{L^{\infty}}\lesssim  \sup_{r\in[0,T]}\norm[\bigg]{\int_{0}^{r}P_{r-s}^{\sigma}[f_{R}(u^{R}_s)-f_{R}(u^{R}_{k_{M}(s)})]ds}_{L^{\infty}}.
\end{align}
In order to obtain the final bound for $e^{M,R}_t$ any $t\in[0,T]$, we observe that 
\begin{align}\label{eq:eMR-decomp}
    e^{M,R}_t=P_{t-k_{M}(t)}^{\sigma}e^{M,R}_{k_{M}(t)} +\int_{k_{M}(t)}^{t}P_{t-s}^{\sigma}[f_R(u^{R}_s)-f_R(u^R_{k_{M}(s)})]ds+\int_{k_{M}(t)}^{t}P_{t-s}^{\sigma}[f_R(u^{R}_{k_M(s)})-f_R(V^{M,R}_{k_{M}(s)})]ds
\end{align} for $t\in[k_{M}(t),T]$.
We first bound $\norm{P_{t-k_{M}(t)}^{\sigma}e^{M,R}_{k_{M}(t)}}_{L^{\infty}}\leqslant \norm{e^{M,R}_{k_{M}(t)}}_{L^{\infty}}$. Using the semigroup estimates \eqref{eq:HK-1} with $d<\sigma<d+2\alpha$ and $\sigma\geq 2\alpha$ together with \eqref{eq:fR-growth} for $m=\lceil \frac{\sigma-d}{2}-\epsilon\rceil$, we find that
\begin{align}\label{eq:a1}
\MoveEqLeft
    \norm[\bigg]{\int_{k_{M}(t)}^{t}P_{t-s}^{\sigma}[f_R(u^{R}_{k_M(s)})-f_R(V^{M,R}_{k_{M}(s)})]ds}_{L^{\infty}}\nonumber
    \\&\lesssim \int_{k_{M}(t)}^t (t-s)^{\frac{\sigma-d-2\alpha-\epsilon}{2\sigma}}\paren[\Big]{\norm{f_R(u^{R}_{k_M(s)})}_{\calC^{\frac{\sigma-d-2\alpha-\epsilon}{2}}}+\norm{f_R(V^{M,R}_{k_{M}(s)})}_{\calC^{\frac{\sigma-d-2\alpha-\epsilon}{2}}}}ds\nonumber
    \\&\lesssim C_R\paren[\Big]{1+\norm{u^{R}}_{C_{T}\calC^{\frac{\sigma-d-\epsilon}{2}}}^m+\norm{V^{M,R}}_{C_{T}\calC^{\frac{\sigma-d-\epsilon}{2}}}^m} \int_{k_{M}(t)}^{t}(t-s)^{\frac{\sigma-d-2\alpha-\epsilon}{2\sigma}}ds\nonumber
    \\&\lesssim M^{-\frac{3\sigma-d-2\alpha-\epsilon}{2\sigma}} C_R\paren[\Big]{1+\norm{u^{R}}_{C_{T}\calC^{\frac{\sigma-d-\epsilon}{2}}}^m+\norm{V^{M,R}}_{C_{T}\calC^{\frac{\sigma-d-\epsilon}{2}}}^m}
\end{align}
and similar we find
\begin{align}\label{eq:a2}
    \norm[\bigg]{\int_{k_{M}(t)}^{t}P_{t-s}^{\sigma}[f_R(u^{R}_s)-f_R(u^R_{k_{M}(s)})]ds}_{L^{\infty}}
    \lesssim M^{-\frac{3\sigma-d-2\alpha-\epsilon}{2\sigma}} C_R\paren[\Big]{1+\norm{u^{R}}_{C_{T}\calC^{\frac{\sigma-d-\epsilon}{2}}}^m}.
\end{align}
By triangle inequality using the bounds \eqref{eq:Gb}, \eqref{eq:a1} and \eqref{eq:a2} for bounding \eqref{eq:eMR-decomp} this yields that
\begin{align*}
  \sup_{t\in[0,T]}\norm{e^{M,R}_{t}}_{L^{\infty}}&\lesssim \sup_{t\in[0,T]}\norm{e^{M,R}_{k_{M}(t)}}_{L^{\infty}} + M^{-\frac{3}{2}+\frac{2\alpha+d}{2\sigma} +\epsilon} C_R\paren[\Big]{1+\norm{u^{R}}_{C_{T}\calC^{\frac{\sigma-d-\epsilon}{2}}}^m+\norm{V^{M,R}}_{C_{T}\calC^{\frac{\sigma-d-\epsilon}{2}}}^m}\\&\lesssim  \sup_{r\in[0,T]}\norm[\bigg]{\int_{0}^{r}P_{r-s}^{\sigma}[f_{R}(u^{R}_s)-f_{R}(u^{R}_{k_{M}(s)})]ds}_{L^{\infty}} \\&\qquad +  M^{-\frac{3}{2}+\frac{2\alpha+d}{2\sigma} +\epsilon} C_R\paren[\Big]{1+\norm{u^{R}}_{C_{T}\calC^{\frac{\sigma-d-\epsilon}{2}}}^m+\norm{V^{M,R}}_{C_{T}\calC^{\frac{\sigma-d-\epsilon}{2}}}^m}.
\end{align*}
Taking the $p$-th moment we obtain the claim due to the a priori estimates \eqref{eq:uR-moment} and \eqref{eq:VRM-moment}, as well as an application of \cref{prop:sewing}. 
\end{proof}
\begin{proof}[Proof of \cref{prop:sewing}]
In order to prove \eqref{eq:sewing}, it suffices (due to Kolmogorov continuity theorem as stated in \cite[Proposition 4.6]{DjGK} and Besov embeddings) to prove that for any $j\geq -1$, $x\in\mathbb{T}$, $0\leq s\leq t\leq T'\leq T$, 
\begin{align}\label{eq:sewing2}
    \E\abs[\bigg]{\int_{s}^{t}\Delta_j P_{T'-r}^{\sigma}[f_{R}(u^{R}_r)-f_{R}(u^{R}_{k_{M}(r)})](x)dr}^p\leq C 2^{-j\epsilon p}M^{p(-\frac{3}{2}+\frac{d+2\alpha}{2\sigma}+2\epsilon)}\abs{t-s}^{\epsilon p}.
\end{align}
To prove \eqref{eq:sewing2}, we apply \cref{lem:ss}. Similar as in the proof of \cite[Proposition 4.5 and Proposition 6.5]{DjGK}, we write $u^{R}=d^{R} + \tilde{O}$, where $\tilde{O}_t:=P_t^{\sigma} u_0 + O_t^{\sigma}$, and define the germ, for fixed $x\in\mathbb{T}^d$ and $j\geq -1$ and $T'\leq T$, as follows, for $0\leq s\leq t\leq T'$, 
\begin{align*}
    A_{s,t}=\int_{s}^{t}\E_s \Delta_j P_{T'-r}^{\sigma} [f_{R}(\E_sd^{R}_r + \tilde O_{r})-f_{R}(\E_s d^{R}_{k_{M}(r)}+\tilde O_{k_{M}(r)})] (x)dr.
\end{align*}
First, we would like to deduce that by the uniqueness of the sewing integral satisfying \eqref{eq:ss1} and \eqref{eq:ss2}, the above germ approximates the desired integral $$\mathcal{A}=\int_{0}^{\cdot}P_{T'-s}^{\sigma}[f_{R}(u^{R}_s)-f_{R}(u^{R}_{k_{M}(s)})]ds.$$ To see this, we have to verify \eqref{eq:ss1} and \eqref{eq:ss2} for $\mathcal{A}_{s,t}-A_{s,t}$. To that aim, we use the intermediate germ 
\begin{align*}
  \tilde{A}_{s,t}=\int_s^t \E_s \Delta_j P_{T'-r}^{\sigma} [f_{R}(d^{R}_r + \tilde O_{r})-f_{R}(d^{R}_{k_{M}(r)}+\tilde O_{k_{M}(r)})] (x)dr
\end{align*}
and verify both bounds for $\mathcal{A}_{s,t}-\tilde{A}_{s,t}$, as well as for $A_{s,t}-\tilde{A}_{s,t}$.
We have that the difference $\mathcal{A}_{s,t}-\tilde{A}_{s,t}$ satisfies \eqref{eq:ss2} with $K_2=0$. Moreover by a trivial estimate, using triangle inequality, semigroup estimates, as well as that the bound on $f_R$ from \eqref{eq:fR-bdd}, we obtain that 
\begin{align*}
    (\E\abs{\mathcal{A}_{s,t}-\tilde{A}_{s,t}}^p)^{1/p}\leqslant (\E\abs{\mathcal{A}_{s,t}}^p)^{1/p}+(\E\abs{\tilde{A}_{s,t}}^p)^{1/p}\lesssim C_R (t-s) (T'-t)^{-\frac{\alpha+\epsilon}{\sigma}} 
\end{align*}
and thus \eqref{eq:ss1} holds true, where we note that $\frac{\alpha+\epsilon}{\sigma}\leq \frac{1}{2}$ for small enough $\epsilon>0$. The verification of \eqref{eq:ss1} for the difference $\tilde{A}_{s,t}-A_{s,t}$ is similar. As for \eqref{eq:ss2}, we first derive regularity bounds on $d_r^{R}-\mathcal{R}_{s,r}$, where $\mathcal{R}_{s,r}:=\E_s d^{R}_r$, that we will need several times in the sequel.
Using \eqref{eq:dR-moment} from \cref{lem:conv-R}, we have that for $q\geq 1$, $\epsilon\in (0,\frac{\sigma-d}{2})$, $0\leq s\leq r\leq T'$,
\begin{align}\label{eq:r1}
    (\E\norm{d_r^{R}-\mathcal{R}_{s,r}}_{L^{\infty}}^q)^{1/q}\lesssim \abs{r-s}^{\frac{3\sigma-2\alpha-d-2\epsilon}{2\sigma}}
\end{align} and, similar as for \cite[equation (6.21)]{DjGK}, we see that
\begin{align}\label{eq:r2}
    (\E\norm{d_{k_{M}(r)}^{R}-\mathcal{R}_{s,k_{M}(r)}}_{L^{\infty}}^q)^{1/q}\lesssim \abs{r-s}^{\frac{3\sigma-2\alpha-d-2\epsilon}{2\sigma}}.
\end{align}
To see \eqref{eq:ss2} for $\tilde{A}_{s,t}- A_{s,t}$, we observe that with \eqref{eq:r1} and \eqref{eq:r2},
\begin{align*}
  (\E\abs{\tilde{A}_{s,t}-A_{s,t}}^p)^{1/p}&\lesssim  C_R (t-s) (T'-t)^{-\frac{\alpha+\epsilon}{\sigma}} \\&\quad\times \sup_{r\in[s,t]}\paren[\big]{(\E\norm{d^{R}_r - \mathcal{R}_{s,r}}_{L^{\infty}}^{p})^{1/p}+\E\norm{d^{R}_{k_{M}(r)} - \mathcal{R}_{s,k_{M}(r)}}_{L^{\infty}}^{p})^{1/p}}
    \\&\lesssim (t-s)^{\frac{5\sigma-2\alpha-d-2\epsilon}{2\sigma}} (T'-t)^{-\frac{\alpha+\epsilon}{\sigma}}.
\end{align*}
Since $3\sigma>2\alpha+d$ as $\sigma\geq 2\alpha$ and $\sigma>d$, together this yields \eqref{eq:ss2} for $\tilde{A}_{s,t}-A_{s,t}$. Thus we conclude that the Riemann sums for the germ $A_{s,t}$ approximate the desired integral $\mathcal{A}$.\\
Now, to apply \cref{lem:ss} to the germ $A_{s,t}$ in order to see \eqref{eq:sewing2}, we are left to prove that for small $\epsilon>0$ and $0\leq s\leq t\leq T'\leq T$, 
\begin{align}\label{eq:g1}
    \E[\abs{A_{s,t}}^p]^{1/p}\lesssim 2^{-j\epsilon}M^{-\frac{3}{2}+\frac{d+2\alpha}{2\sigma}+\frac{3\epsilon}{\sigma}}(T'-t)^{-\frac{1}{2}}(t-s)^{\frac{1}{2}+\frac{\epsilon}{\sigma}} 
\end{align}
and for $s\leq u\leq t$,
\begin{align}\label{eq:g2}
    \E[\abs{\E_s\delta A_{s,u,t}}^p]^{1/p}\lesssim 2^{-j\epsilon}M^{-\frac{3}{2}+\frac{d+2\alpha}{2\sigma}+\frac{3\epsilon}{\sigma}}(T'-t)^{-\frac{1}{2}}(t-s)^{1+\frac{\epsilon}{\sigma}}.
\end{align}
First, we again derive regularity bounds that we use in the sequal. We note that both $\tilde{O}$ and $O$ satisfy the regularity bounds \eqref{eq:ou-reg} and \eqref{eq:s-t-bound} (taking $t=k_{M}(r)-s$ and $s=r-k_{M}(r)$ there) from \cref{prop:ou}, since by the semigroup estimates \eqref{eq:HK-1} for $\hat\lambda,\lambda\in[0,\sigma)$, $\theta\in [0,\sigma-\lambda]$, $s\leq k_{M}(r)\leq r$
\begin{align}
    \norm{(P_{t}^{\sigma}-P_{s}^{\sigma})u_0}_{\calC^{\frac{\sigma-d}{2}-\eps-\hat\lambda}}&\lesssim \abs{t-s}^{\frac{\hat\lambda}{\sigma}}\norm{u_0}_{\calC^{\frac{\sigma-d}{2}-\eps}},\\ 
    \norm{(P_{r-s}^{\sigma}-P_{k_{M}(r)-s}^{\sigma})P^{\sigma}_s u_0}_{\calC^{\frac{\sigma-d}{2}+\eps-\lambda}}&=\norm{(P_{r-k_{M}(r)}^{\sigma}-\operatorname{Id})P^{\sigma}_{k_{M}(r)-s} P^{\sigma}_s u_0}_{\calC^{\frac{\sigma-d}{2}+\eps-\hat\lambda}}\nonumber
    \\&\lesssim (r-k_{M}(r))^{\frac{\lambda+\theta}{\sigma}}(k_{M}(r)-s)^{-\frac{\theta}{\sigma}}\norm{u_0}_{\calC^{\frac{\sigma-d}{2}-\eps}}.
\end{align}
In the following we may thus safely apply the bounds \eqref{eq:ou-reg} and \eqref{eq:s-t-bound} taking $t=k_{M}(r)-s$ and $s=r-k_{M}(r)$ there, for $\tilde O$ in place of $O$.
Using triangle inequality, Jensen's inequality for the conditional expectation, \eqref{eq:dR-moment} and \eqref{eq:ou-reg} taking $\lambda=\sigma-\frac{d}{2}-\alpha-2\epsilon$ there, we have that for $\mathcal{R}_{u,r}=\E_u d^{R}_{r}$,
\begin{align}\label{eq:b1}
    \E[\norm{\mathcal{R}_{s,r}+\tilde O_r-(\mathcal{R}_{s,k_{M}(r)}+\tilde O_{k_{M}(r)})}_{\calC^{\alpha-\sigma/2+\epsilon}}^p]^{1/p}\lesssim (r-k_{M}(r))^{1-\frac{2\alpha+d}{2\sigma} -\frac{2\epsilon}{\sigma}}.
\end{align}
If $k_{M}(r)>s$, we find by triangle inequality and \eqref{eq:dR-moment}, as well as \eqref{eq:s-t-bound} applied for $\theta=\frac{\sigma}{2}-\epsilon$ and $\lambda=\sigma-\alpha-\frac{d}{2} - 2\epsilon$,
\begin{align}\label{eq:b2}
\MoveEqLeft
   \E[\norm{\mathcal{R}_{s,r}+P_{r-s}^{\sigma}\tilde O_s-(\mathcal{R}_{s,k_{M}(r)}+P_{k_{M}(r)-s}^{\sigma}\tilde O_{s})}_{\calC^{\alpha-\sigma/2+\epsilon}}^p]^{1/p}\nonumber\\&\lesssim (r-k_{M}(r))^{\frac{3}{2} -\frac{2\alpha+d}{2\sigma}-\frac{3\epsilon}{\sigma}} (k_{M}(r)-s)^{-\frac{1}{2}+\frac{\epsilon}{\sigma}}.
\end{align}
To prove \eqref{eq:g1}, we distinguish between the cases $\abs{t-s}\leq 3M^{-1}$ and $\abs{t-s}> 3M^{-1}$. 
In the first case, we bound with \eqref{eq:fR-negLip} applied for $s=\frac{\sigma-2\alpha}{2}$ using that by \cref{ass:1}, $2\alpha> d$,
\begin{align*}
\MoveEqLeft
    \abs{A_{s,t}}\lesssim 2^{-j\epsilon}\norm[\bigg]{\int_{s}^{t}\E_s P_{T'-r}^{\sigma} [f_{R}(\E_sd^{R}_r + \tilde O_{r})-f_{R}(\E_s d^{R}_{k_{M}(r)}+\tilde O_{k_{M}(r)})] dr}_{\calC^{\epsilon}}
    \\&\lesssim 2^{-j\epsilon} \E_s\int_{s}^{t}(T'-r)^{-\frac{1}{2}}\norm{f_R(\mathcal{R}_{s,r} + \tilde O_{r})-f_{R}(\mathcal{R}_{s,k_{M}(r)}+\tilde O_{k_{M}(r)})}_{\calC^{-\frac{\sigma}{2}+\epsilon}}dr
    \\&\lesssim 2^{-j\epsilon} C_{R}\E_s\paren[\bigg]{\paren[\big]{1+\norm{d^{R}}_{C_T\calC^{\frac{\sigma}{2}-\alpha}}^m+\norm{\tilde{O}}_{C_T\calC^{\frac{\sigma}{2}-\alpha}}^m}\times\\&\hspace{3cm}\int_{s}^{t}(T'-r)^{-\frac{1}{2}}\norm{\mathcal{R}_{s,r} + \tilde O_{r}-(\mathcal{R}_{s,k_{M}(r)}+\tilde O_{k_{M}(r)})}_{\calC^{-\frac{\sigma}{2}+\alpha+\epsilon}}dr},
\end{align*}
where $m=\lceil \frac{\sigma-2\alpha}{2}\rceil$.
Using \eqref{eq:b1}, as well as the moment bounds for $\tilde{O}$ from \eqref{eq:ou-reg} and $d^{R}$ from \eqref{eq:dR-moment}, and tower property for the conditional expectation, we obtain since $\abs{t-s}\leq 3M^{-1}$,
\begin{align*}
\MoveEqLeft
\E[\abs{A_{s,t}}^p]^{1/p}
    \\&\lesssim  2^{-j\epsilon} (T'-t)^{-\frac{1}{2}} (t-s) (r-k_{M}(r))^{1-\frac{2\alpha+d}{2\sigma} -2\epsilon/\sigma} 
   \\& = 2^{-j\epsilon} (T'-t)^{-\frac{1}{2}} (t-s)^{\frac{1}{2}+\frac{\epsilon}{\sigma}} (t-s)^{\frac{1}{2}-\frac{\epsilon}{\sigma}} (r-k_{M}(r))^{1-\frac{2\alpha+d-4\eps}{2\sigma}} 
    \\&\lesssim 2^{-j\epsilon} (T'-t)^{-\frac{1}{2}} (t-s)^{\frac{1}{2}+\frac{\epsilon}{\sigma}} M^{-\frac{3}{2} +\frac{2\alpha+d}{2\sigma}+\frac{3\epsilon}{\sigma}}.
\end{align*}
In the case of $\abs{t-s}> 3M^{-1}$, we define $t'$ as the second smallest grid point bigger than $s$, such that $k_{M}(r)-u\geq (r-u)/2$ for all $r\geq t'$. We decompose the integral $\int_{s}^{t}$ into $\int_{s}^{t'}$ for which the same argument as above applies, since $\abs{t'-s}\leq 2M^{-1}$, and $\int_{t'}^{t}$ for which we argue as follows. To compute the conditional expectation, we use that the nonlinearity satisfies $f_R(u)=\mathcal{D}_{\alpha} G_R(u)$ for a Nemytskii operator $G_R$. Indeed, for a continuous $G:\R\to\R$ and a $\F_s$-measurable random variable $X$ and $s\leq r$, we may compute for any $x\in\T$, almost surely
\begin{align*}
    \E_s[G(X+O_r^{\sigma})] (x) &= \E_s [G(X+O_r^{\sigma}(x))]
    = (P_{Q(r-s)}^{\R}G) (X+P_{r-s}^{\sigma} O_s^{\sigma})(x),
\end{align*} where $Q(r-s)=\E(\int_{s}^{r}p^{\sigma}(r-v,x-y)\xi(dv,dy))^2$ is the variance, which is constant in $x$ and $P^{\R}$ denotes the heat semigroup on $\R$. By continuity in $x$, we furthermore obtain that, almost surely, as processes,
\begin{align*}
    \E_s[G(X+O_r^{\sigma})]=(P_{Q(r-s)}^{\R}G) (X+P_{r-s}^{\sigma} O_s^{\sigma}).
\end{align*}
Exchanging the conditional expectation with the derivative and applying the above for $G=G_R$, we compute
\begin{align*}
\MoveEqLeft
\int_{t'}^{t}\E_s \Delta_j P_{T'-r}^{\sigma} [f_{R}(\E_sd^{R}_r + \tilde O_{r})-f_{R}(\E_s d^{R}_{k_{M}(r)}+\tilde O_{k_{M}(r)})] dr
\\&= \int_{t'}^{t} \Delta_j \mathcal{D}_{\alpha} P_{T'-r}^{\sigma} \E_s [G_{R}(\E_sd^{R}_r + \tilde O_{r})-G_{R}(\E_s d^{R}_{k_{M}(r)}+\tilde O_{k_{M}(r)})] dr
\\& = \int_{t'}^{t} \Delta_j \mathcal{D}_{\alpha} P_{T'-r}^{\sigma} [(P^{\R}_{Q(r-s)}G_{R})(\mathcal{R}_{s,r} + P^{\sigma}_{r-s}\tilde O_{s})-(P_{Q(k_{M}(r)-s)}^{\R}G_{R})(\mathcal{R}_{s,k_{M}(r)}+P^{\sigma}_{k_{M}(r)-s}\tilde O_{s})] dr
\\& =  \int_{t'}^{t} \Delta_j \mathcal{D}_{\alpha} P^{\sigma}_{T'-r} [(P^{\R}_{Q(r-s)}G_{R})(\mathcal{R}_{s,r} + P^{\sigma}_{r-s}\tilde O_{s})-(P^{\R}_{Q(r-s)}G_{R})(\mathcal{R}_{s,k_{M}(r)}+P^{\sigma}_{k_{M}(r)-s}\tilde O_{s})] dr
\\&\qquad + \int_{t'}^{t} \Delta_j \mathcal{D}_{\alpha} P^{\sigma}_{T'-r} [(P^{\R}_{Q(r-s)}G_{R}-P_{Q(k_{M}(r)-s)}^{\R}G_{R})(\mathcal{R}_{s,r} + P^{\sigma}_{r-s}\tilde O_{s})]dr
\\&=: I_{1} + I_{2}
\end{align*}
We bound $\E\abs{I_{1}(x)}^p$ similar as in \cite[equation (4.23)]{DjGK}.
Indeed by \eqref{eq:composition}, we see that for $u\in\calC^{s}$ and $s\in (0,\lceil\frac{\sigma-d}{2}\rceil)\setminus\N$,
\begin{align}\label{eq:GR-growth2}
    \norm{(P^{\R}_{Q(r-s)}G_R)' (u)}_{\calC^{s}}\lesssim C_R (1+\norm{u}_{\calC^{s}}^{\lceil s\rceil}).
\end{align}
Using the latter bound for $s=\frac{\sigma-2\alpha}{2}$ and $m=\lceil s\rceil$ with $2\alpha>d$ by Assumption \ref{ass:1} and \eqref{eq:b2}, we proceed to estimate $I_1$ as follows
\begin{align*}
(&\E\abs{I_{1}(x)}^p)^{1/p}\nonumber
        \\&=\paren[\bigg]{\E\abs[\bigg]{\int_{t'}^{t}\Delta_j \mathcal{D}_{\alpha} P^{\sigma}_{T' - r}[(P^{\R}_{Q(r-u)}G_R)(\mathcal{R}_{s,r}+P^{\sigma}_{r-s} \tilde O_{s})-(P^{\R}_{Q(r-s)}G_R)(\mathcal{R}_{s,k_{M}(r)}+P^{\sigma}_{k_{M}(r)-u} \tilde O_{u})](x)dr}^{p}}^{1/p}\nonumber
        \\&\lesssim 2^{-j\epsilon}\paren[\bigg]{\E\norm[\bigg]{\int_{t'}^{t} \mathcal{D}_{\alpha} P^{\sigma}_{T' - r}\bigg[\paren[\bigg]{\int_{0}^{t}(P^{\R}_{Q(r-s)}G_R)'(\lambda (\mathcal{R}_{s,r}+P^{\sigma}_{r-s}\tilde O_{s})+(1-\lambda)(\mathcal{R}_{s,k_{M}(r)}+P^{\sigma}_{k_{M}(r)-s}\tilde O_s))d\lambda} \nonumber
        \\&\hspace{5cm}\times(\mathcal{R}_{s,r}+P^{\sigma}_{r-s} \tilde O_{s}-\mathcal{R}_{s,k_{M}(r)}-P^{\sigma}_{k_{M}(r)-s} \tilde O_{s})\bigg]dr}_{\calC^{\epsilon}}^{p}}^{1/p}\nonumber
        \\&\lesssim 2^{-j\epsilon}\int_{t'}^{t}\paren[\bigg]{(T' - r)^{-\frac{1}{2}}\nonumber
        \\&\hspace{1cm}\times\paren[\bigg]{\E\norm[\bigg]{\int_{0}^{1}(P^{\R}_{Q(r-s)}G_R)'(\lambda (\mathcal{R}_{s,r}+P^{\sigma}_{r-s}\tilde O_{s})+(1-\lambda)(\mathcal{R}_{s,k_{M}(r)}+P^{\sigma}_{k_{M}(r)-s}\tilde O_s))d\lambda}_{\calC^{\frac{\sigma-2\alpha}{2}}}^{2p}}^{1/2p}}\nonumber
        \\&\hspace{1cm}\times\paren[\big]{\paren{\E\norm{\mathcal{R}_{s,r}+P^{\sigma}_{r-s} \tilde O_{s}-\mathcal{R}_{s,k_{M}(r)}-P^{\sigma}_{k_{M}(r)-s} \tilde O_{s}}_{\calC^{\frac{2\alpha-\sigma}{2}+\epsilon}}^{2p}}^{1/2p}}dr\nonumber
        \\&\lesssim 2^{-j\epsilon} C_R(\E(1+\sup_{0\leq s'\leq r'\leq T'}\norm{\mathcal{R}_{s',r'}+\tilde{O}_{r'}}_{\calC^{\frac{\sigma-2\alpha}{2}}}^{m})^{2p})^{1/2p})\\&\qquad\times
        \int_{t'}^{t}(T' - r)^{-\frac{1}{2}} (r-k_{M}(r))^{\frac{3}{2}-\frac{2\alpha+d}{2\sigma}-\frac{3\epsilon}{\sigma}}(k_{M}(r)-s)^{-\frac{1}{2}+\frac{\epsilon}{\sigma}}dr\nonumber
        \\&\lesssim 2^{-j\epsilon}M^{-\frac{3}{2}+\frac{2\alpha+d}{2\sigma}+\frac{3\epsilon}{\sigma}}
        \int_{t'}^{t}(T' - r)^{-\frac{1}{2}} (r-s)^{-\frac{1}{2}+\frac{\epsilon}{\sigma}}dr\nonumber
        \\&\lesssim 2^{-j\epsilon}M^{-\frac{3}{2}+\frac{2\alpha+d}{2\sigma}+\frac{3\epsilon}{\sigma}} (T'-t)^{-\frac{1}{2}}(t-s)^{\frac{1}{2}+\frac{\epsilon}{\sigma}},
    \end{align*}
    where we used that $k_{M}(r)-s\geq (r-s)/2$ for $r\in [t',t]$ in the penultimate inequality. 
    For $I_{2}$, we use the following estimate on $Q$ (cf. \cite[equation (4.25)]{DjGK}, but using $\norm{p^{\sigma}_s}_{L^{\infty}(\T^d)}\lesssim s^{-\frac{d}{\sigma}}$), where $s\leq l\leq r$ and $\epsilon\in (0,\frac{\sigma}{2})$,
    \begin{align}\label{eq:Q-est}
        Q(r-s)-Q(l-s)
        \lesssim (r-l)^{\frac{3}{2}-\frac{d+\epsilon}{\sigma}}(l-s)^{-\frac{1}{2}+\frac{\epsilon}{\sigma}}.
    \end{align}
     The latter bound together with the heat kernel estimate \eqref{eq:HK} yields
     \begin{align*}
         \norm{(P^{\R}_{Q(r-s)}G_R)-(P^{\R}_{Q(k_{M}(r)-s)}G_R)}_{C^{0}_{b}}&\lesssim (Q(r-s)-Q(k_{M}(r)-s))\norm{G_{R}}_{C^{2}_{b}}\\&\lesssim (k_{M}(r)-s)^{-\frac{1}{2}+\frac{\epsilon}{\sigma}}(r-k_{M}(r))^{\frac{3}{2}-\frac{d+\epsilon}{\sigma}}.
     \end{align*} 
    Hence for $I_2$ we obtain the estimate,
    \begin{align*}
    (&\E\abs{I_{2}(x)}^p)^{1/p}\nonumber
        \\&=\paren[\bigg]{\E\abs[\bigg]{\int_{t'}^{t}\Delta_j \mathcal{D}_{\alpha} P^{\sigma}_{T' - r}\paren[\Big]{P^{\R}_{Q(r-s)}G_R-P^{\R}_{Q(k_{M}(r)-s)}G_R}(\mathcal{R}_{s,k_{M}(r)}+P^{\sigma}_{k_{M}(r)-s}\tilde O_{s})(x)dr}^p}^{1/p}\nonumber
        \\&\lesssim 2^{-j\epsilon}\paren[\bigg]{\int_{t'}^{t}(T' - r)^{-\frac{\alpha+\epsilon}{\sigma}}\norm{P^{\R}_{Q(r-s)}G_R-P^{\R}_{Q(k_{M}(r)-s)}G_R}_{L^{\infty}}dr}
        \\&\lesssim  2^{-j\epsilon}(T'-t)^{-\frac{\alpha+\epsilon}{\sigma}}C_R\int_{t'}^{t}(k_{M}(r)-s)^{-\frac{1}{2}+\frac{\epsilon}{\sigma}}(r-k_{M}(r))^{\frac{3}{2}-\frac{d+\epsilon}{\sigma}}dr\nonumber
        \\&\lesssim 2^{-j\epsilon}M^{-\frac{3}{2}+\frac{d}{\sigma}+\frac{\epsilon}{\sigma}} (T'-t)^{-\frac{\alpha+\epsilon}{\sigma}}(t-s)^{\frac{1}{2}+\frac{\epsilon}{\sigma}},
    \end{align*}
    using again that $k_{M}(r)-s\geq (r-s)/2$ for $r\in [t',t]$ in the last inequality. Since $\sigma> 2\alpha$, we note that $\frac{\alpha+\epsilon}{\sigma}\leq\frac{1}{2}$ for small enough $\epsilon>0$. Since $d\leq 2\alpha$, the rate for $I_2$ furthermore dominates the one of $I_1$. Thus the estimate for $I_1$ and $I_2$ together yield \eqref{eq:g1}. 

To prove \eqref{eq:g2}, we write 
\begin{align*}
    \E_s \delta A_{sut} &= \Delta_j \int_{u}^{t} \E_s\mathcal{D}_{\alpha}\E_u P^{\sigma}_{T' - r}\paren[\Big]{G_R(\mathcal{R}_{s,r} + \tilde O_{r})-G_R(\mathcal{R}_{s,k_{M}(r)} + \tilde O_{r})
    \\&\hspace{3cm} - (G_R(\mathcal{R}_{u,r} + \tilde O_{r})-G_R(\mathcal{R}_{u,k_{M}(r)} + \tilde O_{r}))} (x) dr
\end{align*}
We can rewrite the difference inside integral in the following form:
\begin{align*}
\MoveEqLeft
    G_R(\cR_{s,r}+\tilde O_r)-G_R(\cR_{s,k_{M}(r)}+\tilde O_{k_{M}(r)})
    -(G_R(\cR_{u,r}+\tilde O_r)-G_R(\cR_{u,k_{M}(r)}+\tilde O_{k_{M}(r)})) \\&= (\cR_{s,r}-\cR_{u,r})\int_{0}^{1}G'_{R}(\lambda (\cR_{s,r}+\tilde O_{r})+(1-\lambda)(\cR_{u,r}+\tilde O_r))d\lambda
    \\&\quad +(\cR_{s,k_{M}(r)}-\cR_{u,k_{M}(r)})
    \\&\hspace{3cm}\times\squeeze[0.5]{\int_{0}^{1}G'_{R}(\lambda (\cR_{s,k_{M}(r)}+\tilde O_{k_{M}(r)})+(1-\lambda)(\cR_{u,k_{M}(r)}+\tilde O_{k_{M}(r)}))d\lambda}
    \\&= (\cR_{s,r}-\cR_{u,r})\int_{0}^{1}\paren[\bigg]{G'_{R}(\lambda (\cR_{s,r}+\tilde O_{r})+(1-\lambda)(\cR_{u,r}+\tilde O_r))
    \\&\hspace{3.5cm}-\squeeze[0.5]{G'_{R}(\lambda (\cR_{s,k_{M}(r)}+\tilde O_{k_{M}(r)})+(1-\lambda) (\cR_{u,k_{M}(r)}+\tilde O_{k_{M}(r)}))}}d\lambda
    \\&\quad + (\cR_{s,k_{M}(r)}-\cR_{s,r}+\cR_{u,r}-\cR_{u,k_{M}(r)})
    \\&\hspace{3cm}\times\squeeze[0.5]{\int_{0}^{1}G'_{R}(\lambda (\cR_{s,k_{M}(r)}+\tilde{O}_{k_{M}(r)})+(1-\lambda)(\cR_{u,k_{M}(r)}+\tilde O_{k_{M}(r)}))d\lambda}
    \\&=:J_1+J_2.
\end{align*}
From \eqref{eq:r1} and triangle inequality, we deduce that for $s\leq u\leq r$,
\begin{align}\label{eq:r11a}
    (\E\norm{\mathcal{R}_{s,r}-\mathcal{R}_{u,r}}_{L^{\infty}}^q)^{1/q}\lesssim \abs{r-s}^{\frac{3\sigma-2\alpha-d-2\epsilon}{2\sigma}}
\end{align} and similar, using again \eqref{eq:dR-moment} and triangle inequality, we deduce
\begin{align}\label{eq:r11b}
    (\E\norm{\mathcal{R}_{s,r}-\mathcal{R}_{u,r}}_{\calC^{\sigma/2-\alpha}}^q)^{1/q}\lesssim \abs{r-s}^{\frac{2\sigma-d-4\epsilon}{2\sigma}}.
\end{align} 
Furthermore, by \eqref{eq:r2} with triangle inequality, we find that
\begin{align}\label{eq:r22}
    (\E\norm{\mathcal{R}_{s,k_{M}(r)}-\mathcal{R}_{u,k_{M}(r)}}_{L^{\infty}}^q)^{1/q}\lesssim_q \abs{r-s}^{\frac{3\sigma-2\alpha-d-2\epsilon}{2\sigma}}.
\end{align}
Moreover by Jensen's inequality for the conditional expectation and \cref{lem:remainder} we have that
\begin{align}\label{eq:r33}
   (\E\norm{\mathcal{R}_{s,r}-\mathcal{R}_{s,k_{M}(r)}}_{L^{\infty}}^q)^{1/q}\lesssim_q \abs{r-k_{M}(r)}^{\frac{3\sigma-2\alpha-d-2\epsilon}{2\sigma}}.
\end{align}
Since $\norm{G'_R(u)}_{L^{\infty}}\leq C_R$ and using once \eqref{eq:r11a} and \eqref{eq:r22} and once \eqref{eq:r33}, we obtain for $J_2$ the following estimate 
\begin{align*}
\MoveEqLeft
    \paren[\bigg]{\E\abs[\bigg]{\Delta_j \int_{u}^{t}\E_{s}[\mathcal{D}_{\alpha} P^{\sigma}_{T' - r}\E_{u}[J_2]](x)dr}^{p}}^{1/p}\\&\lesssim
    2^{-j\epsilon}(T'-t)^{-\frac{\alpha-\epsilon}{\sigma}}C_R
    \paren[\bigg]{\E\paren[\bigg]{\int_{u}^{t}\norm{\cR_{s,k_{M}(r)}-\cR_{s,r}+\cR_{u,r}-\cR_{u,k_{M}(r)}}_{L^{\infty}}dr}^{2p}}^{1/2p}
    \\&\lesssim 2^{-j\epsilon}(T'-t)^{-\frac{\alpha-\epsilon}{\sigma}}
    \int_{u}^{t}\min((r-k_{M}(r))^{\frac{3\sigma-2\alpha-d-2\epsilon}{2\sigma}}, (r-s)^{\frac{3\sigma-2\alpha-d-2\epsilon}{2\sigma}})dr
    \\&\lesssim 2^{-j\epsilon}(T'-t)^{-\frac{\alpha-\epsilon}{\sigma}}
    (t-s)^{1+\frac{\epsilon}{\sigma}}M^{\frac{2\alpha-3\sigma+d+4\epsilon}{2\sigma}},
\end{align*} 
where the last bound follows by the interpolation $\min(a,b)\leq a^{1-\theta}b^{\theta}$, $a,b\geq 0$. 

To estimate $J_{1}$ we distinguish again the cases $\abs{t-u}\leq 3M^{-1}$ and $\abs{t-u}>3M^{-1}$ like for $A_{s,t}$.
In the case of $\abs{t-u}\leq 3M^{-1}$, the estimate for $J_1$ is similar to the estimate for $A_{s,t}$ above in the case of $\abs{t-s}\leq 3M^{-1}$ but for $G'_R$ in place of $G_R$ using that $G_{R}\in C^{3}_{b}$ by Assumption \ref{ass:G}, since $3\leq \lceil\frac{\sigma-d}{2}\rceil +2$. However in $J_1$ we have the additional factor $\cR_{s,r}-\cR_{u,r}$, which did not appear in $A_{s,t}$ and for which we use the estimate \eqref{eq:r11b}. Together we obtain that if $\abs{t-u}\leq 3M^{-1}$,
\begin{align*}
\paren[\bigg]{\E\abs[\bigg]{\Delta_j \int_{u}^{t}\E_{s}[\mathcal{D}_{\alpha} P^{\sigma}_{T' - r}\E_{u}[J_1]](x)dr}^{p}}^{1/p}\lesssim     2^{-j\epsilon}M^{-\frac{3}{2}+\frac{2\alpha+d}{2\sigma}+\frac{3\epsilon}{\sigma}} (T'-t)^{-\frac{1}{2}}(t-u)^{\frac{3}{2}-\frac{d}{2\sigma}-\frac{2\epsilon}{\sigma}}.
\end{align*}
If $\abs{t-u}>3M^{-1}$, we introduce $t'$ as before, in particular $k_{M}(r)\geq u\geq s$.
The integral from $u$ to $t'$ is treated as in the case $|t-u|\leq 3M^{-1}$. For the integral from $t'$ to $t$ we have that, using that the factor $\cR_{s,r}-\cR_{u,r}$ is $\F_{u}$-measurable,
\begin{align*}
    &\paren[\bigg]{\E\abs[\bigg]{\Delta_j \int_{t'}^{t}\E_{s}[\mathcal{D}_{\alpha} P^{\sigma}_{T' - r}\E_{u}[J_1]](x)dr}^{p}}^{1/p}\\&\lesssim 2^{-j\epsilon}(T'-t)^{-\frac{1}{2}} \paren[\bigg]{\E\paren[\bigg]{\int_{t'}^{t}\norm{\cR_{s,r}-\cR_{u,r}}_{\calC^{\frac{\sigma-2\alpha}{2}}}\\&\quad\times\int_{0}^{1}\norm[\Big]{(P^{\R}_{Q(r-u)}G_{R}')(\lambda \cR_{s,r}+(1-\lambda)\cR_{u,r}+P^{\sigma}_{r-u}\tilde O_u)
    \\&\qquad\quad\squeeze[0.5]{-(P^{\R}_{Q(k_{M}(r)-u)}G_{R}')(\lambda \cR_{s,k_{M}(r)}+(1-\lambda)\cR_{u,k_{M}(r)}+P^{\sigma}_{k_{M}(r)-u}\tilde O_u)}}_{\calC^{\frac{2\alpha-\sigma}{2}+\epsilon}} d\lambda dr}^{p}}^{1/p}.
\end{align*}
Using similar steps as for estimating the time integral from $t'$ to $t$ above, but now applied for $G_{R}'\in C^{\lceil \frac{\sigma-d}{2}\rceil +1}_{b}$ instead of $G_{R}$, to bound
\begin{align*}
\MoveEqLeft
    \paren[\bigg]{\E\int_{t'}^{t}\norm{(P^{\R}_{Q(r-u)}G_{R}')(\lambda \cR_{s,r}+(1-\lambda)\cR_{u,r}+P^{\sigma}_{r-u}\tilde O_u)
    \\&\hspace{0.6cm}\squeeze[0.5]{-(P^{\R}_{Q(k_{M}(r)-u)}G_{R}')(\lambda \cR_{s,k_{M}(r)}+(1-\lambda)\cR_{u,k_{M}(r)}+P^{\sigma}_{k_{M}(r)-u}\tilde O_u)}}_{\calC^{\frac{\sigma-2\alpha}{2}}}^{2p} dr}^{1/2p}
\end{align*}
for any fixed $\lambda\in[0,1]$, as well as \eqref{eq:r11b} and \eqref{eq:b2}, we arrive at
\begin{align*}
\MoveEqLeft
    \paren[\bigg]{\E\abs[\bigg]{\Delta_j \int_{t'}^{t}\E_{s}[\mathcal{D}_{\alpha} P^{\sigma}_{T' - r}\E_{u}[J_1]](x)dr}^{p}}^{1/p}\\&\lesssim
    2^{-j\epsilon}(t-s)^{\frac{2\sigma-d-\eps}{2\sigma}}(T'-t)^{-\frac{1}{2}}
    \int_{t'}^{t}(r-k_{M}(r))^{\frac{3}{2}-\frac{2\alpha+d}{2\sigma}-\frac{3\epsilon}{\sigma}}(k_{M}(r)-u)^{-\frac{1}{2}+\frac{\epsilon}{\sigma}}dr
    \\&\lesssim 2^{-j\epsilon}(t-s)^{\frac{4\sigma-d-2\eps}{2\sigma}}(T'-t)^{-\frac{1}{2}}M^{-\frac{3}{2}+\frac{2\alpha+d}{2\sigma}+\frac{3\epsilon}{\sigma}} .
\end{align*}
The estimates for $J_1$ and $J_2$ together yield \eqref{eq:g2}, which concludes the proof.
\end{proof}

\section{Examples}\label{sec:ex}
In this section, we prove that \cref{ass:1} is satisfied in the cases of interest: Burgers equation \eqref{ex:Burgers} and Cahn-Hilliard equation \eqref{ex:CH}.

\begin{proposition}\label{prop:wp}
Let for any $\epsilon>0$, $u_{0}\in\calC^{\frac{\sigma-d}{2}-\epsilon}$ and assume \eqref{ex:Burgers} or \eqref{ex:CH}. Then there exists a unique global mild solution $u$ to \eqref{eq:SPDE}.
Moreover for any $\lambda\in [0,\frac{\sigma}{2})$ and any small enough $\epsilon>0$, it follows that $u\in C_{T}^{\frac{\lambda}{\sigma}}\calC^{\frac{\sigma-d}{2}-\lambda-\epsilon}$ almost surely. 
\end{proposition}
\begin{proof}
The existence and uniqueness of a global mild solution to \eqref{eq:SPDE} with $u\in C_T L^2$ almost surely is known by \cite{DaPrato} in the case of \eqref{ex:Burgers}. 
In the case of \eqref{ex:CH}, we apply \cite[Theorem 1.3]{CW01} for $q$ arbitrarily large and obtain a global mild solution $u$ with $u\in L^{\infty}([0,T],L^{q}(\mathbb{T}^d))$ almost surely. Notice that \cite[Theorem 1.3]{CW01} is formulated for Neumann boundary conditions, but a complete analogue proof holds for periodic boundary conditions.
Given the mild solution $u$ with $u\in C_{T}L^{2}$, respectively $u\in L^{\infty}_TL^{q}$ almost surely, we employ semigroup estimates and embeddings to improve the spatial and temporal regularity. Recall that $f(u)=\mathcal{D}_{\alpha} G(u)$. 
Since $G$ is a polynomial of the general form $G(x)=\sum_{k=0}^{m}c_k x^k$ for $m\in\N$ and constants $c_k\in\R$, by \eqref{eq:derivative} and \eqref{eq:product} for $s>0$ and $q\geq m$, we obtain that for $u\in B_{q,\infty}^{s}$,
\begin{align}\label{eq:f-bound1}
    \norm{f(u)}_{B^{s-\alpha}_{q/m,\infty}}\lesssim\sum_{k=0}^{m}\norm{u^k}_{B^{s}_{q/m,\infty}}\lesssim 1+\norm{u}_{B_{q,\infty}^{s}}^{m}
\end{align}
and for $u\in L^q$
\begin{align}\label{eq:f-bound2}
    \norm{f(u)}_{B^{-\alpha}_{q/m,\infty}}\lesssim\sum_{k=0}^{m}\norm{u^k}_{L^{q/m}}\lesssim 1+\norm{u}_{L^{q}}^{m}.
\end{align}
Using the embedding $ B^{\frac{\sigma-d}{2}+\frac{d(m-1)}{q}-\epsilon}_{q/m,\infty}\hookrightarrow B^{\frac{\sigma-d}{2}-\epsilon}_{q,\infty}$, as well as the semigroup estimates together with $1-\frac{\sigma-d+2\alpha-2\epsilon}{2\sigma}-\frac{d(m-1)}{q\sigma}>0$ (which is true for $\epsilon>0$, since for \eqref{ex:Burgers} we have $q=2, m=2, d=1, \sigma=2\alpha=2$ and for \eqref{ex:CH} we have $\sigma=2\alpha=4$, $m=3$ and choose $q> 8$) and \eqref{eq:f-bound2}, we have that for any $t\in[0,T]$ and  $\epsilon>0$ small enough,
\begin{align}\label{eq:u-apriori1}
    \norm{u_t}_{B^{\frac{\sigma-d}{2}-\epsilon}_{q,\infty}}
    &\lesssim\norm{P_{t}^{\sigma}u_{0}}_{B^{\frac{\sigma-d}{2}-\eps}_{q,\infty}}+\int_{0}^{t}\norm{P_{t-s}^{\sigma}f(u_s)}_{B^{\frac{\sigma-d}{2}+\frac{d(m-1)}{q}-\epsilon}_{q/m,\infty}}ds + \norm{O^{\sigma}}_{C_T\calC^{\frac{\sigma-d}{2}-\epsilon}}\nonumber
    \\&\lesssim \norm{u_{0}}_{\calC^{\frac{\sigma-d}{2}-\eps}}+\int_{0}^{t}(t-s)^{-\frac{\sigma-d+2\alpha-2\epsilon}{2\sigma}-\frac{d(m-1)}{q\sigma}}\norm{f(u_s)}_{B^{-\alpha}_{q/m,\infty}}ds + \norm{O^{\sigma}}_{C_T\calC^{(\sigma-d)/2-\epsilon}}\nonumber
    \\&\lesssim \norm{u_{0}}_{\calC^{\frac{\sigma-d}{2}-\eps}}+T^{1-\frac{\sigma-d+2\alpha-2\epsilon}{2\sigma}-\frac{d(m-1)}{q\sigma}}(1+\norm{u}_{L^{\infty}_T L^{q}}^{m}) +\norm{O^{\sigma}}_{C_T\calC^{(\sigma-d)/2-\epsilon}}.
\end{align}
Moreover, using the embedding $B^{\frac{\sigma-d}{2}+\frac{dm}{q}-\epsilon}_{q/m,\infty}\hookrightarrow \calC^{\frac{\sigma-d}{2}-\epsilon}$ and semigroup estimates together with $1-\frac{\alpha}{\sigma}-\frac{dm}{q\sigma}+\frac{\epsilon}{\sigma}>0$ (which is true since $\epsilon>0$ and for \eqref{ex:Burgers} we have $\sigma=2\alpha=2, d=1, m=2, q=2$ and for \eqref{ex:CH} we have $\sigma=4, \alpha=2, m=3$ and choose $q>d$), we obtain together with \eqref{eq:f-bound1} for any $t\in[0,T]$  and  small enough $\epsilon>0$,
\begin{align}\label{eq:reg-bound}
    \norm{u_t}_{\calC^{\frac{\sigma-d}{2}-2\epsilon}}&\lesssim \norm{u_{0}}_{\calC^{\frac{\sigma-d}{2}-\eps}}+\int_{0}^{t}(t-s)^{-\frac{\alpha}{\sigma}-\frac{dm}{q\sigma}+\frac{\epsilon}{\sigma}}\norm{f(u_s)}_{B^{\frac{\sigma-d}{2}-\alpha-\epsilon}_{q/m,\infty}}ds + \norm{O^{\sigma}}_{C_T\calC^{\frac{\sigma-d}{2}-\epsilon}}\nonumber
    \\&\lesssim_T \norm{u_{0}}_{\calC^{\frac{\sigma-d}{2}-\eps}}+T^{1-\frac{\alpha}{\sigma}-\frac{dm}{q\sigma}+\frac{\epsilon}{\sigma}}(1+ \norm{u}_{L^{\infty}_T B^{\frac{\sigma-d}{2}-\epsilon}_{q,\infty}}^{m}) +\norm{O^{\sigma}}_{C_T\calC^{\frac{\sigma-d}{2}-\epsilon}}.
\end{align}
Plugging \eqref{eq:u-apriori1} into the right-hand side of \eqref{eq:reg-bound} and using \cref{prop:ou}, we find that for any $\epsilon>0$ small enough, $u\in L^{\infty}_T\calC^{\frac{\sigma-d}{2}-\epsilon}$ almost surely.
To see the time regularity, we write $u_t=P_{t}^{\sigma}u_0 + d_t + O^{\sigma}_t$, $t\in[0,T]$, where we know that by \cref{prop:ou} that for $\lambda\in [0,\frac{\sigma}{2})$, $O^{\sigma}\in C_T^{\frac{\lambda}{\sigma}}\calC^{\frac{\sigma-d}{2}-\lambda-\epsilon}$ almost surely. By the semigroup estimates we have that for $\lambda\in [0,\sigma]$, $$\norm{P_{t}^{\sigma}u_{0}-P_{s}^{\sigma}u_{0}}_{\calC^{\frac{\sigma-d}{2} -\lambda-\eps}}\lesssim \abs{t-s}^{\frac{\lambda}{\sigma}}\norm{u_0}_{\calC^{\frac{\sigma-d}{2}-\eps}}.$$ 
Thus it is left to bound $d$. We have that by the semigroup estimates for $\lambda\in [0,\sigma)$, using that $\sigma\geq 2\alpha\geq\alpha$, and \eqref{eq:f-bound1},
\begin{align*}
\MoveEqLeft
    \norm{d_t - d_s}_{\calC^{\frac{\sigma-d}{2}-\lambda-2\epsilon}}\nonumber
    \\&\lesssim\norm{d_t - d_s}_{\calC^{\frac{3\sigma-d}{2}-\alpha-\lambda-2\epsilon}}\nonumber
    \\&\leq \int_{s}^{t}\norm{P_{t-r}^{\sigma} f(u_r)}_{\calC^{\frac{3\sigma-d}{2}-\alpha-\lambda-2\epsilon}}dr + \int_0^s \norm{P_{s-r}^{\sigma}(\operatorname{Id}-P_{t-s}^{\sigma})f(u_r)}_{\calC^{\frac{3\sigma-d}{2}-\alpha-\lambda-2\epsilon}}dr\nonumber
    \\&\lesssim \int_{s}^{t} (t-r)^{-1+\frac{\lambda+\eps}{\sigma}}\norm{f(u_r)}_{\calC^{\frac{\sigma-d-2\eps}{2}-\alpha}}dr \nonumber
    + \int_0^s (t-s)^{\frac{\lambda}{\sigma}}(s-r)^{-1+\frac{\epsilon}{2\sigma}}\norm{f(u_r)}_{\calC^{\frac{\sigma-d-2\eps}{2}-\alpha}} dr \nonumber
    \\&\lesssim (t-s)^{\frac{\lambda}{\sigma}}(1+\norm{u}_{L^{\infty}_T\calC^{\frac{\sigma-d}{2}-\epsilon}}^{m}).
\end{align*}
Since by the above for any $\epsilon>0$ small enough, $u\in L^{\infty}_T\calC^{\frac{\sigma-d}{2}-\eps} $ almost surely, we deduce that for any $\lambda\in [0,\sigma)$, $\epsilon>0$ small, $d\in C_T^{\frac{\lambda}{\sigma}}\calC^{\frac{3\sigma-d}{2}-\alpha-\lambda-\epsilon}$ almost surely and for any $\lambda\in [0,\frac{\sigma}{2})$, $\epsilon>0$ small, $u\in C_T^{\frac{\lambda}{\sigma}}\calC^{\frac{\sigma-d}{2}-\alpha-\lambda-\epsilon} $ almost surely.
\end{proof}

\newcommand{\etalchar}[1]{$^{#1}$}


\begin{thebibliography}{ABLM24}
\expandafter\ifx\csname url\endcsname\relax
  \def\url#1{\texttt{#1}}\fi
\expandafter\ifx\csname urlprefix\endcsname\relax\def\urlprefix{URL }\fi
\expandafter\ifx\csname href\endcsname\relax
  \def\href#1#2{#2}\fi
\expandafter\ifx\csname burlalt\endcsname\relax
  \def\burlalt#1#2{\href{#2}{\texttt{#1}}}\fi

\bibitem[ABLM24]{ABLM}
\textsc{S.~Athreya}, \textsc{O.~Butkovsky}, \textsc{K.~Lê}, and
  \textsc{L.~Mytnik}.
\newblock Well-posedness of stochastic heat equation with distributional drift
  and skew stochastic heat equation.
\newblock \emph{Communications on Pure and Applied Mathematics} \textbf{77},
  no.~5, (2024), 2708--2777.
\newblock
  \burlalt{doi:https://doi.org/10.1002/cpa.22157}{http://dx.doi.org/https://doi.org/10.1002/cpa.22157}.

\bibitem[AG06]{Alabert2006}
\textsc{A.~Alabert} and \textsc{I.~Gyöngy}.
\newblock \emph{On Numerical Approximation of Stochastic {Burgers'} Equation},
  1--15.
\newblock Springer Berlin Heidelberg, 2006.

\bibitem[BCD11]{Bahouri2011}
\textsc{H.~Bahouri}, \textsc{J.-Y. Chemin}, and \textsc{R.~Danchin}.
\newblock \emph{Fourier Analysis and Nonlinear Partial Differential Equations}.
\newblock Grundlehren der mathematischen Wissenschaften. Springer Berlin
  Heidelberg, 2011.
\newblock
  \burlalt{doi:10.1007/978-3-642-16830-7}{http://dx.doi.org/10.1007/978-3-642-16830-7}.

\bibitem[BCM26]{brehier2024weak}
\textsc{C.-E. Br{\'e}hier}, \textsc{S.~Cox}, and \textsc{A.~Millet}.
\newblock Weak error estimates of {G}alerkin approximations for the stochastic
  {Burgers} equation driven by additive trace-class noise.
\newblock \emph{Potential Analysis} \textbf{64}, no.~2, (2026), 35.

\bibitem[BCW26]{brehier2022weak}
\textsc{C.-E. Br{\'e}hier}, \textsc{J.~Cui}, and \textsc{X.~Wang}.
\newblock Weak error estimates of fully-discrete schemes for the stochastic
  {Cahn}--{Hilliard} equation.
\newblock \emph{Communications in Mathematical Sciences} \textbf{24}, no.~2,
  (2026), 521--564.

\bibitem[BFG22]{BFG}
\textsc{C.~Bellingeri}, \textsc{P.~K. Friz}, and \textsc{M.~Gerencsér}.
\newblock Singular paths spaces and applications.
\newblock \emph{Stochastic Analysis and Applications} \textbf{40}, no.~6,
  (2022), 1126--1149.
\newblock
  \burlalt{doi:10.1080/07362994.2021.1988641}{http://dx.doi.org/10.1080/07362994.2021.1988641}.

\bibitem[BGJK20]{BGJK}
\textsc{S.~Becker}, \textsc{B.~Gess}, \textsc{A.~Jentzen}, and
  \textsc{P.~Kloeden}.
\newblock Lower and upper bounds for strong approximation errors for numerical
  approximations of stochastic heat equations.
\newblock \emph{BIT Numerical Mathematics} \textbf{60}, (2020), 1057--–1073.
\newblock
  \burlalt{doi:10.1007/s10543-020-00807-2}{http://dx.doi.org/10.1007/s10543-020-00807-2}.

\bibitem[BGJK23]{becker2023strong}
\textsc{S.~Becker}, \textsc{B.~Gess}, \textsc{A.~Jentzen}, and \textsc{P.~E.
  Kloeden}.
\newblock Strong convergence rates for explicit space-time discrete numerical
  approximations of stochastic allen-cahn equations.
\newblock \emph{Stochastics and Partial Differential Equations: Analysis and
  Computations} \textbf{11}, no.~1, (2023), 211--268.
\newblock
  \burlalt{doi:10.1007/s40072-021-00226-6}{http://dx.doi.org/10.1007/s40072-021-00226-6}.

\bibitem[BHJ{\etalchar{+}}19]{beccari2019strong}
\textsc{M.~Beccari}, \textsc{M.~Hutzenthaler}, \textsc{A.~Jentzen},
  \textsc{R.~Kurniawan}, \textsc{F.~Lindner}, and \textsc{D.~Salimova}.
\newblock Strong and weak divergence of exponential and linear-implicit euler
  approximations for stochastic partial differential equations with
  superlinearly growing nonlinearities.
\newblock \emph{arXiv preprint arXiv:1903.06066} (2019).

\bibitem[BJ13]{BJ13}
\textsc{D.~Bl\"{o}mker} and \textsc{A.~Jentzen}.
\newblock Galerkin approximations for the stochastic {Burgers} equation.
\newblock \emph{SIAM Journal on Numerical Analysis} \textbf{51}, no.~1, (2013),
  694--715.
\newblock \burlalt{doi:10.1137/110845756}{http://dx.doi.org/10.1137/110845756}.

\bibitem[BLR25]{blomker2025numerical}
\textsc{D.~Bl{\"o}mker}, \textsc{C.~Ling}, and \textsc{J.~Rimmele}.
\newblock Numerical approximation of {Cahn}--{Hilliard} type nonlinear {SPDEs}
  with additive space-time white noise.
\newblock \emph{arXiv preprint arXiv:2501.18240} (2025).

\bibitem[BM25]{BanasMukam}
\textsc{L.~Ba{\v{n}}as} and \textsc{J.~D. Mukam}.
\newblock Numerical approximation of the stochastic {Cahn}-–{Hilliard}
  equation with space–time white noise near the sharp-interface limit.
\newblock \emph{IMA Journal of Numerical Analysis} (2025).
\newblock
  \burlalt{doi:10.1093/imanum/draf064}{http://dx.doi.org/10.1093/imanum/draf064}.

\bibitem[Br{\'e}24]{Brehier2024}
\textsc{C.-E. Br{\'e}hier}.
\newblock Total variation error bounds for the accelerated exponential {Euler}
  scheme approximation of parabolic semilinear {SPDEs}.
\newblock \emph{SIAM Journal on Numerical Analysis} \textbf{62}, no.~3, (2024),
  1171--1190.
\newblock
  \burlalt{doi:10.1137/22M152596X}{http://dx.doi.org/10.1137/22M152596X}.

\bibitem[BYZ23]{BanasYangZhu2023}
\textsc{{\v{L}}.~Ba{\v{n}}as}, \textsc{H.~Yang}, and \textsc{R.~Zhu}.
\newblock Sharp interface limit of stochastic {Cahn--Hilliard} equation with
  singular noise.
\newblock \emph{Potential Analysis} \textbf{59}, (2023), 497--518.
\newblock
  \burlalt{doi:10.1007/s11118-021-09976-3}{http://dx.doi.org/10.1007/s11118-021-09976-3}.

\bibitem[CCHS22]{ChenCuiHongSheng2022}
\textsc{C.~Chen}, \textsc{J.~Cui}, \textsc{J.~Hong}, and \textsc{D.~Sheng}.
\newblock Accelerated exponential {Euler} scheme for stochastic heat equation:
  convergence rate of the density.
\newblock \emph{IMA Journal of Numerical Analysis} \textbf{43}, no.~2, (2022),
  1181--1220.
\newblock
  \burlalt{doi:10.1093/imanum/drac011}{http://dx.doi.org/10.1093/imanum/drac011}.

\bibitem[CHS20]{cui2020}
\textsc{J.~Cui}, \textsc{J.~Hong}, and \textsc{L.~Sun}.
\newblock Numerical analysis of a full discretization for stochastic
  {Cahn}--{Hilliard} equation driven by additive noise.
\newblock \emph{arXiv preprint arXiv:1812.06289} (2020).

\bibitem[Coo70]{Cook1970}
\textsc{H.~E. Cook}.
\newblock Brownian motion in spinodal decomposition.
\newblock \emph{Acta Metallurgica} \textbf{18}, no.~3, (1970), 297--306.
\newblock
  \burlalt{doi:10.1016/0001-6160(70)90144-6}{http://dx.doi.org/10.1016/0001-6160(70)90144-6}.

\bibitem[CQW23]{cai2023strong}
\textsc{M.~Cai}, \textsc{R.~Qi}, and \textsc{X.~Wang}.
\newblock Strong convergence rates of an explicit scheme for stochastic
  {Cahn}--{Hilliard} equation with additive noise.
\newblock \emph{BIT Numerical Mathematics} \textbf{63}, no.~4, (2023), 43.

\bibitem[CW01]{CW01}
\textsc{C.~Cardon-Weber}.
\newblock Cahn-hilliard stochastic equation: existence of the solution and of
  its density.
\newblock \emph{Bernoulli} \textbf{7}, no.~5, (2001), 777--816.

\bibitem[DG01]{Davie-Gaines}
\textsc{A.~M. Davie} and \textsc{J.~G. Gaines}.
\newblock Convergence of numerical schemes for the solution of parabolic
  stochastic partial differential equations.
\newblock \emph{Mathematics of Computation} \textbf{70}, no. 233, (2001),
  121--135.
\newblock
  \burlalt{doi:10.1090/s0025-5718-00-01224-2}{http://dx.doi.org/10.1090/s0025-5718-00-01224-2}.

\bibitem[DGK26]{DjGK}
\textsc{A.~Djurdjevac}, \textsc{M.~Gerencs{\'e}r}, and \textsc{H.~Kremp}.
\newblock {Higher order approximation of nonlinear SPDEs with additive
  space-time white noise}.
\newblock \emph{The Annals of Probability} \textbf{54}, no.~6, (2026), 2788 --
  2826.
\newblock
  \burlalt{doi:10.1214/25-AOP1818}{http://dx.doi.org/10.1214/25-AOP1818}.

\bibitem[DGL23]{DGL}
\textsc{K.~Dareiotis}, \textsc{M.~Gerencsér}, and \textsc{K.~Lê}.
\newblock {Quantifying a convergence theorem of Gy\"ongy and Krylov}.
\newblock \emph{The Annals of Applied Probability} \textbf{33}, no.~3, (2023),
  2291--2323.
\newblock
  \burlalt{doi:10.1214/22-aap1867}{http://dx.doi.org/10.1214/22-aap1867}.

\bibitem[DPD96]{DaPratoDebussche}
\textsc{G.~Da~Prato} and \textsc{A.~Debussche}.
\newblock Stochastic {Cahn}-{Hilliard} equation.
\newblock \emph{Nonlinear Anal., Theory Methods Appl.} \textbf{26}, no.~2,
  (1996), 241--263.
\newblock
  \burlalt{doi:10.1016/0362-546X(94)00277-O}{http://dx.doi.org/10.1016/0362-546X(94)00277-O}.

\bibitem[DPDT94]{DaPrato}
\textsc{G.~Da~Prato}, \textsc{A.~Debussche}, and \textsc{R.~Temam}.
\newblock Stochastic {Burgers'} equation.
\newblock \emph{NoDEA} \textbf{1}, (1994), 389–402.
\newblock
  \burlalt{doi:https://doi.org/10.1007/BF01194987}{http://dx.doi.org/https://doi.org/10.1007/BF01194987}.

\bibitem[DZ07]{DebusscheZambotti2007}
\textsc{A.~Debussche} and \textsc{L.~Zambotti}.
\newblock Conservative stochastic {Cahn--Hilliard} equation with reflection.
\newblock \emph{The Annals of Probability} \textbf{35}, no.~5, (2007),
  1706--1739.
\newblock
  \burlalt{doi:10.1214/009117906000000773}{http://dx.doi.org/10.1214/009117906000000773}.

\bibitem[EM91]{ElezovicMikelic1991}
\textsc{N.~Elezovi{\'c}} and \textsc{A.~Mikeli{\'c}}.
\newblock On the stochastic {Cahn--Hilliard} equation.
\newblock \emph{Nonlinear Analysis: Theory, Methods \& Applications}
  \textbf{16}, no.~12, (1991), 1169--1200.
\newblock
  \burlalt{doi:10.1016/0362-546X(91)90204-E}{http://dx.doi.org/10.1016/0362-546X(91)90204-E}.

\bibitem[GIP15]{GIP}
\textsc{M.~Gubinelli}, \textsc{P.~Imkeller}, and \textsc{N.~Perkowski}.
\newblock {Paracontrolled distributions and singular {PDEs}}.
\newblock \emph{Forum of Mathematics, Pi} \textbf{3}, (2015), e6.
\newblock
  \burlalt{doi:10.1017/fmp.2015.2}{http://dx.doi.org/10.1017/fmp.2015.2}.

\bibitem[Gy{\"o}99]{Gy}
\textsc{I.~Gy{\"o}ngy}.
\newblock {Lattice Approximations for Stochastic Quasi-Linear Parabolic Partial
  Differential Equations driven by Space-Time White Noise II}.
\newblock \emph{Potential Analysis} \textbf{11}, no.~1, (1999), 1--37.
\newblock
  \burlalt{doi:10.1023/a:1008699504438}{http://dx.doi.org/10.1023/a:1008699504438}.

\bibitem[Hai23]{hairer2023}
\textsc{M.~Hairer}.
\newblock An introduction to stochastic {PDEs}.
\newblock \emph{arXiv:0907.4178} (2023).

\bibitem[Hen81]{Henry}
\textsc{D.~Henry}.
\newblock \emph{Geometric theory of semilinear parabolic equations}, vol. 840.
\newblock Springer Berlin Heidelberg New York, 1981.

\bibitem[HJ20]{HutzenthalerJentzenzbMATH07206753}
\textsc{M.~Hutzenthaler} and \textsc{A.~Jentzen}.
\newblock On a perturbation theory and on strong convergence rates for
  stochastic ordinary and partial differential equations with nonglobally
  monotone coefficients.
\newblock \emph{Ann. Probab.} \textbf{48}, no.~1, (2020), 53--93.
\newblock
  \burlalt{doi:10.1214/19-AOP1345}{http://dx.doi.org/10.1214/19-AOP1345}.

\bibitem[HJLP19]{hutzenthaler2019strong}
\textsc{M.~Hutzenthaler}, \textsc{A.~Jentzen}, \textsc{F.~Lindner}, and
  \textsc{P.~Pu{\v{s}}nik}.
\newblock Strong convergence rates on the whole probability space for
  space-time discrete numerical approximation schemes for stochastic burgers
  equations.
\newblock \emph{arXiv preprint arXiv:1911.01870} (2019).

\bibitem[HJS19]{HutzenthalerJentzenSalimova2018}
\textsc{M.~Hutzenthaler}, \textsc{A.~Jentzen}, and \textsc{D.~Salimova}.
\newblock Strong convergence of full-discrete nonlinearity-truncated
  accelerated exponential {Euler}-type approximations for stochastic
  {Kuramoto}-{Sivashinsky} equations.
\newblock \emph{Commun. Math. Sci.} \textbf{16}, no.~6, (2019), 1489--1529.
\newblock
  \burlalt{doi:10.4310/CMS.2018.v16.n6.a2}{http://dx.doi.org/10.4310/CMS.2018.v16.n6.a2}.

\bibitem[HL22]{hutzenthaler2022strong}
\textsc{M.~Hutzenthaler} and \textsc{R.~Link}.
\newblock Strong convergence rates for full-discrete approximations of
  stochastic {Burgers} equations with multiplicative noise.
\newblock \emph{arXiv preprint arXiv:2210.17536} (2022).

\bibitem[HM12]{HairerMaas12_zbMATH06067453}
\textsc{M.~Hairer} and \textsc{J.~Maas}.
\newblock A spatial version of the {It{\=o}}-{Stratonovich} correction.
\newblock \emph{Ann. Probab.} \textbf{40}, no.~4, (2012), 1675--1714.
\newblock \burlalt{doi:10.1214/11-AOP662}{http://dx.doi.org/10.1214/11-AOP662}.

\bibitem[HM16]{HairerMatetski_zbMATH06600872}
\textsc{M.~Hairer} and \textsc{K.~Matetski}.
\newblock Optimal rate of convergence for stochastic {Burgers}-type equations.
\newblock \emph{Stoch. Partial Differ. Equ., Anal. Comput.} \textbf{4}, no.~2,
  (2016), 402--437.
\newblock
  \burlalt{doi:10.1007/s40072-015-0067-5}{http://dx.doi.org/10.1007/s40072-015-0067-5}.

\bibitem[HV11]{HairerVoss_zbMATH06031653}
\textsc{M.~Hairer} and \textsc{J.~Voss}.
\newblock Approximations to the stochastic {Burgers} equation.
\newblock \emph{J. Nonlinear Sci.} \textbf{21}, no.~6, (2011), 897--920.
\newblock
  \burlalt{doi:10.1007/s00332-011-9104-3}{http://dx.doi.org/10.1007/s00332-011-9104-3}.

\bibitem[Jen11]{Jentzen12}
\textsc{A.~Jentzen}.
\newblock Higher order pathwise numerical approximations of {SPDEs} with
  additive noise.
\newblock \emph{{SIAM} Journal on Numerical Analysis} \textbf{49}, no.~2,
  (2011), 642--667.
\newblock \burlalt{doi:10.1137/080740714}{http://dx.doi.org/10.1137/080740714}.

\bibitem[JK08]{Jentzen-Kloden}
\textsc{A.~Jentzen} and \textsc{P.~E. Kloeden}.
\newblock Overcoming the order barrier in the numerical approximation of
  stochastic partial differential equations with additive
  space{\textendash}time noise.
\newblock \emph{Proceedings of the Royal Society A: Mathematical, Physical and
  Engineering Sciences} \textbf{465}, no. 2102, (2008), 649--667.
\newblock
  \burlalt{doi:10.1098/rspa.2008.0325}{http://dx.doi.org/10.1098/rspa.2008.0325}.

\bibitem[JSW19]{JSW19}
\textsc{A.~Jentzen}, \textsc{D.~Salimova}, and \textsc{T.~Welti}.
\newblock Strong convergence for explicit space–time discrete numerical
  approximation methods for stochastic {Burgers} equations.
\newblock \emph{Journal of Mathematical Analysis and Applications}
  \textbf{469}, no.~2, (2019), 661--704.
\newblock
  \burlalt{doi:https://doi.org/10.1016/j.jmaa.2018.09.032}{http://dx.doi.org/https://doi.org/10.1016/j.jmaa.2018.09.032}.

\bibitem[KH12]{MinooHosseini_zbMATH06107055}
\textsc{M.~Kamrani} and \textsc{S.~M. Hosseini}.
\newblock Spectral collocation method for stochastic {Burgers} equation driven
  by additive noise.
\newblock \emph{Math. Comput. Simul.} \textbf{82}, no.~9, (2012), 1630--1644.
\newblock
  \burlalt{doi:10.1016/j.matcom.2012.03.007}{http://dx.doi.org/10.1016/j.matcom.2012.03.007}.

\bibitem[L{\^e}20]{Khoa}
\textsc{K.~L{\^e}}.
\newblock A stochastic sewing lemma and applications.
\newblock \emph{Electronic Journal of Probability} \textbf{25}, (2020), 1--55.
\newblock
  \burlalt{doi:doi:10.1214/20-EJP442}{http://dx.doi.org/doi:10.1214/20-EJP442}.

\bibitem[MS20]{MazzonettoSalimova20}
\textsc{S.~Mazzonetto} and \textsc{D.~Salimova}.
\newblock Existence, uniqueness, and numerical approximations for stochastic
  {Burgers} equations.
\newblock \emph{Stochastic Anal. Appl.} \textbf{38}, no.~4, (2020), 623--646.
\newblock
  \burlalt{doi:10.1080/07362994.2019.1709503}{http://dx.doi.org/10.1080/07362994.2019.1709503}.

\bibitem[Pri01]{Printems}
\textsc{J.~Printems}.
\newblock On the discretization in time of parabolic stochastic partial
  differential equations.
\newblock \emph{ESAIM: Mathematical Modelling and Numerical Analysis}
  \textbf{35}, no.~6, (2001), 1055–1078.
\newblock
  \burlalt{doi:10.1051/m2an:2001148}{http://dx.doi.org/10.1051/m2an:2001148}.

\bibitem[QCW22]{QiCaiWang2022}
\textsc{R.~Qi}, \textsc{M.~Cai}, and \textsc{X.~Wang}.
\newblock Strong convergence rates of a fully discrete scheme for the
  {Cahn--Hilliard--Cook} equation.
\newblock \emph{arXiv preprint arXiv:2204.01630} (2022).
\newblock
  \burlalt{doi:10.48550/arXiv.2204.01630}{http://dx.doi.org/10.48550/arXiv.2204.01630}.

\bibitem[Wan20]{Wang}
\textsc{X.~Wang}.
\newblock An efficient explicit full-discrete scheme for strong approximation
  of stochastic allen{\textendash}cahn equation.
\newblock \emph{Stochastic Processes and their Applications} \textbf{130},
  no.~10, (2020), 6271--6299.
\newblock
  \burlalt{doi:10.1016/j.spa.2020.05.011}{http://dx.doi.org/10.1016/j.spa.2020.05.011}.

\bibitem[WC25]{WangCao}
\textsc{Y.~Wang} and \textsc{W.~Cao}.
\newblock Strong convergence of a fully discrete scheme for stochastic
  {Burgers} equation with fractional-type noise.
\newblock \emph{Adv. Comput. Math.} \textbf{51}, no.~2, (2025), 32.
\newblock Id/No 15.
\newblock
  \burlalt{doi:10.1007/s10444-025-10227-x}{http://dx.doi.org/10.1007/s10444-025-10227-x}.

\bibitem[WCC24]{wang2024strong}
\textsc{Y.~Wang}, \textsc{W.~Cao}, and \textsc{Y.~Cao}.
\newblock Strong convergence of an explicit full-discrete scheme for stochastic
  {Burgers}--{Huxley} equation.
\newblock \emph{arXiv preprint arXiv:2408.00947} (2024).

\end{thebibliography}
\end{document}